\documentclass{amsart}
\usepackage[all]{xy}
\usepackage[mathscr]{eucal}
\usepackage{graphicx}
\CompileMatrices

\newtheorem{theorem}{Theorem}
\newtheorem{lemma}[theorem]{Lemma}
\newtheorem{proposition}[theorem]{Proposition}
\newtheorem{corollary}[theorem]{Corollary}

\theoremstyle{definition}
\newtheorem{definition}[theorem]{Definition}
\newtheorem{example}[theorem]{Example}

\theoremstyle{remark}
\newtheorem{remark}[theorem]{Remark}

\numberwithin{equation}{section}

\begin{document}

\newcommand{\Diff}{\operatorname{Diff}}
\newcommand{\Homeo}{\operatorname{Homeo}}
\newcommand{\Hom}{\operatorname{Hom}}
\newcommand{\Exp}{\operatorname{Exp}}
\newcommand{\Orb}{\operatorname{\textup{Orb}}}
\newcommand{\Stwo}{\mbox{$\displaystyle S^2$}}
\newcommand{\Sn}{\mbox{$\displaystyle S^n$}}
\newcommand{\supp}{\operatorname{supp}}
\newcommand{\intr}{\operatorname{int}}
\newcommand{\kernel}{\operatorname{ker}} \newcommand{\A}{\mathbb{A}}
\newcommand{\B}{\mathbb{B}} \newcommand{\C}{\mathbb{C}}
\newcommand{\D}{\mathbb{D}} \newcommand{\E}{\mathbb{E}}
\newcommand{\F}{\mathbb{F}} \newcommand{\G}{\mathbb{G}}
\newcommand{\Hh}{\mathbb{H}} \newcommand{\I}{\mathbb{I}}
\newcommand{\J}{\mathbb{J}} \newcommand{\K}{\mathbb{K}}
\newcommand{\Ll}{\mathbb{L}} \newcommand{\M}{\mathbb{M}}
\newcommand{\N}{\mathbb{N}} \newcommand{\Oo}{\mathbb{O}}
\newcommand{\Pp}{\mathbb{P}} \newcommand{\Q}{\mathbb{Q}}
\newcommand{\R}{\mathbb{R}} \newcommand{\Ss}{\mathbb{S}}
\newcommand{\T}{\mathbb{T}} \newcommand{\U}{\mathbb{U}}
\newcommand{\V}{\mathbb{V}} \newcommand{\W}{\mathbb{W}}
\newcommand{\X}{\mathbb{X}} \newcommand{\Y}{\mathbb{Y}}
\newcommand{\Z}{\mathbb{Z}} \newcommand{\kk}{\mathbb{k}}
\newcommand{\orbify}[1]{\ensuremath{\mathcal{#1}}}
\newcommand{\Orbdiff}{\ensuremath{\Diff_{\Orb}}}
\newcommand{\RedOrbDiff}{\ensuremath{\Diff_{\textup{red}}}}
\newcommand{\Frechet}{Fr\'{e}chet\ }
\newcommand{\Frechetnospace}{Fr\'{e}chet}

\title[A Manifold Structure for Orbifold Diffeomorphism Groups]{A
  Manifold Structure for the Group of Orbifold Diffeomorphisms of a
  Smooth Orbifold}

\author{Joseph E. Borzellino}
\address{Department of Mathematics, California Polytechnic State
  University, 1 Grand Avenue, San Luis Obispo, California 93407}
\email{jborzell@calpoly.edu}
\thanks{We wish to thank the referee for comments and suggestions that lead to improvements in the final version of this manuscript.}

\author{Victor Brunsden} \address{Department of Mathematics and
  Statistics, Penn State Altoona, 3000 Ivyside Park, Altoona,
  Pennsylvania 16601} \email{vwb2@psu.edu}

\subjclass[2000]{Primary 57S05, 22F50, 54H99; Secondary 22E65}

\keywords{Orbifolds, Diffeomorphism
  Groups, Topological Transformation Groups, Homeomorphism Groups}
\begin{abstract}
  For a compact, smooth $C^r$ orbifold (without boundary), we show
  that the topological structure of the orbifold diffeomorphism group
  is a Banach manifold for $1\le r<\infty$ and a \Frechet manifold if
  $r=\infty$. In each case, the local model is the separable Banach
  (\Frechetnospace) space of $C^r (C^\infty,\textup{resp.})$ orbisections of the
  tangent orbibundle.
\end{abstract}

\maketitle

\section{Introduction}\label{IntroSection}

Our goal in this paper is to determine the topological structure of
the orbifold diffeomorphism group of a smooth compact orbifold.  It is
well known that in the case of a closed smooth $C^r$ manifold, the
group of $C^r$  diffeomorphisms ($1\le r\le\infty$) is a smooth manifold
whose local model is $\mathscr{D}^r(M)$, the space of $C^r$ tangent
vector fields on $M$.  See, for example \cite{MR98h:22024} or
\cite{MR58:31210}.  $\mathscr{D}^r(M)$ is a separable Banach space for
$1\le r < \infty$ and a separable \Frechet space for $r=\infty$. One
might naively think that the orbifold diffeomorphism group is itself
an infinite dimensional orbifold, but one only need remember that the
orbifold diffeomorphism group is a (topological) group and hence must
be homogeneous. As such, it cannot be a non-trivial orbifold. In fact,
in the case of a smooth compact orbifold, the structure of the
orbifold diffeomorphism group holds no surprises, and we have the
following

\begin{theorem}\label{MainTheorem}
  Let $r\ge 1$ and let \orbify{O} be a compact, smooth $C^r$ orbifold
  (without boundary). Denote by $\Diff^r_{\Orb}(\orbify{O})$ the group
  of $C^r$ orbifold diffeomorphisms equipped with the $C^r$ topology.
  Then $\Diff^r_{\Orb}(\orbify{O})$ is a manifold modeled on the
  topological vector space $\mathscr{D}^r_{\Orb}(\orbify{O})$ of $C^r$
  orbisections of the tangent orbibundle equipped with the $C^r$
  topology. This separable vector space is a Banach space if $1\le r <
  \infty$ and is a \Frechet space if $r = \infty$.
\end{theorem}

This particular result was first conjectured with a plausibility argument in \cite{MR2003g:58013}.
Here, we provide a complete proof using techniques in the spirit of the classical result for the manifold case. 
There are many competing and useful notions of smooth orbifold map in the literature.
In \cite{MR2003g:58013}, the statement of theorem~\ref{MainTheorem} referred to {\em unreduced}
orbifold diffeomorphisms. The main result of \cite{MR2003146} concerned the {\em reduced} orbifold diffeomorphism group
$\Diff^r_{\textup{red}}(\orbify{O})$.  It is possible to recover the
topology of $\Diff^r_{\textup{red}}(\orbify{O})$ as a quotient of
$\Diff^r_{\Orb}(\orbify{O})$. In fact, we have the following structure
theorem for $\Diff^r_{\textup{red}}(\orbify{O})$ as a corollary of
theorem~\ref{MainTheorem}.

\begin{corollary}\label{MainCorollary}
  Let $r\ge 1$ and let \orbify{O} be a compact, smooth $C^r$ orbifold
  (without boundary).  Let $\mathscr{ID}=
  \{f\in\Diff^r_{\Orb}(\orbify{O})\ |\ f(x) = x $ for all $
  x\in\orbify{O}\}$. That is, $\mathscr{ID}$ is the set comprised of
  all lifts of the identity map. Then $ |\mathscr{ID}| < \infty $ and
  there is a short exact sequence of groups
  \begin{equation*}
    1\to \mathscr{ID}\to \Diff^r_{\Orb}(\orbify{O})\to
    \Diff^r_{\textup{red}}(\orbify{O})\to 1.
  \end{equation*}
  Thus, $
  \Diff^r_{\textup{red}}(\orbify{O})\cong\Diff^r_{\Orb}(\orbify{O})/\mathscr{ID}$
  is a Banach manifold if $r < \infty $ and a \Frechet manifold if $ r
  = \infty $.
\end{corollary}

\begin{remark} Using methods detailed in \cite{MR1471480}, it will follow that these diffeomorphism groups have the structure of smooth manifolds. Furthermore, composition and inversion in these groups will be continuous, and in the 
$r=\infty$ case, both $\Orbdiff^\infty(\orbify{O})$ and $\RedOrbDiff^\infty(\orbify{O})$ will be convenient \Frechet Lie groups. Details will appear in a future revision to the preprint \cite{BB_Stratified} on the topological structure of the set of smooth mappings between orbifolds $\orbify{O}$ and $\orbify{P}$.
\end{remark}

The next few sections of the paper will define and describe the
notions that appear in the statement of theorem~\ref{MainTheorem} and
corollary~\ref{MainCorollary}. In particular, in
section~\ref{OrbifoldSection}, we define the notion of smooth orbifold
and its natural stratification. We also define the notion of product
orbifold and suborbifold and give some examples. In
section~\ref{OrbifoldMapSection}, we define the notion of orbifold
map.  Section~\ref{CrTopologySection} defines the (strong) $C^r$
topology on maps between smooth orbifolds.  In
section~\ref{TangentObibundleSection}, we define the tangent
orbibundle and its orbisections. The space of orbisections provide the local
model for the orbifold diffeomorphism group. In
section~\ref{SmoothRiemannianStructureSection}, we look at smooth
Riemannian structures and define a smooth Riemannian exponential map.
Finally, we prove theorem~\ref{MainTheorem} and
corollary~\ref{MainCorollary} in
section~\ref{ProofOfMainTheoremSection}.

It should be noted that we have chosen not to use the language of Lie
groupoids and Morita equivalence in our description of orbifolds and
their maps, but rather we have chosen a more ``classical" approach. The reason for this choice is that a treatment using groupoids, in
our opinion, would not add clarity to the exposition or enhance our
results. In fact, we believe that much of the useful geometric and
topological intuition becomes obscured. A reader interested in the groupoid approach to orbifolds and its utility should consult the recent monograph \cite{AdemLeidaRuan} and the references therein, especially the article \cite{MR1950948}. 

We should also note that our definition of orbifold is modeled on the definition in Thurston
\cite{Thurston78}. The orbifolds that concern us here are referred to as {\em classical effective orbifolds} in \cite{AdemLeidaRuan}. While our notion of orbifold map is more general than that given in \cite{AdemLeidaRuan}, our notion of {\em reduced} orbifold map and {\em reduced} orbifold diffeomorphism agrees with its definitions~1.3 and 1.4.

\section{Orbifolds}\label{OrbifoldSection}

In this section, we review the (classical) definition of smooth orbifold and related constructions.

\begin{definition}\label{orbifold}
  An $n$-dimensional (topological) \emph{orbifold} $\orbify{O}$,
  consists of a paracompact, Hausdorff topological space
  $X_\orbify{O}$ called the \emph{underlying space}, with the
  following local structure.  For each $x \in X_\orbify{O}$ and
  neighborhood $U$ of $x$, there is a neighborhood $U_x \subset U$, an
  open set $\tilde U_x \cong \R^n$, a finite group $\Gamma_x$ acting
  continuously and effectively on $\tilde U_x$ which fixes $0\in\tilde
  U_x$, and a homeomorphism $\phi_x:\tilde U_x/\Gamma_x \to U_x$ with
  $\phi_x(0)=x$.  These actions are subject to the condition that for
  a neighborhood $U_z\subset U_x$ with corresponding $\tilde U_z \cong
  \R^n$, group $\Gamma_z$ and homeomorphism $\phi_z:\tilde
  U_z/\Gamma_z \to U_z$, there is an embedding $\tilde\psi_{zx}:\tilde
  U_z \to \tilde U_x$ and an injective homomorphism
  $\theta_{zx}:\Gamma_z \to \Gamma_x$ so that $\tilde\psi_{zx}$ is
  equivariant with respect to $\theta_{zx}$ (that is, for
  $\gamma\in\Gamma_z, \tilde\psi_{zx}(\gamma\cdot \tilde
  y)=\theta_{zx}(\gamma)\cdot\tilde\psi_{zx}(\tilde y)$ for all
  $\tilde y\in\tilde U_z$), such that the following diagram commutes:
  \begin{equation*}
    \xymatrix{{\tilde U_z}\ar[rr]^{\tilde\psi_{zx}}\ar[d]&&{\tilde U_x}\ar[d]\\
      {\tilde U_z/\Gamma_z}\ar[rr]^>>>>>>>>>>{\psi_{zx}=\tilde\psi_{zx}/\Gamma_z}\ar[dd]^{\phi_z}&&{\tilde U_x/\theta_{zx}(\Gamma_z)\ar[d]}\\
      &&{\tilde U_x/\Gamma_x}\ar[d]^{\phi_x}\\
      {U_z}\ar[rr]^{\subset}&&{U_x} }
  \end{equation*}
\end{definition}

\begin{remark}
  Note that if $\delta\in\Gamma_x$ then
  $\overline\psi_{zx}=\delta\cdot\tilde\psi_{zx}$ is also an embedding
  of $\tilde U_z$ into $\tilde U_x$. It is equivariant relative to the
  injective homomorphism
  $\overline\theta_{zx}(\gamma)=\delta\cdot\theta_{zx}(\gamma)\cdot\delta^{-1}$.
  Thus, we regard $\tilde{\psi}_{zx}$ as being defined only up to
  composition with elements of $\Gamma_x$, and $\theta_{zx}$ defined
  only up to conjugation by elements of $\Gamma_x$.  In general, it is
  not true that
  $\tilde{\psi}_{zx}=\tilde{\psi}_{yx}\circ\tilde{\psi}_{zy}$ when
  $U_z\subset U_y\subset U_x$, but there should be an element
  $\delta\in\Gamma_x$ such that
  $\delta\cdot\tilde{\psi}_{zx}=\tilde{\psi}_{yx}\circ\tilde{\psi}_{zy}$
  and $\delta\cdot
  \theta_{zx}(\gamma)\cdot\delta^{-1}=\theta_{yx}\circ\theta_{zy}(\gamma)$.
  Also, the covering $\{U_x\}$ of $X_{\orbify{O}}$ is not an intrinsic
  part of the orbifold structure. We regard two coverings to give the
  same orbifold structure if they can be combined to give a larger
  covering still satisfying the definitions.
\end{remark}

\begin{definition}\label{locallysmooth}
  We say that an $n$-dimensional orbifold $\orbify{O}$ is
  \emph{locally smooth} if the action of $\Gamma_x$ on $\tilde
  U_x\cong\R^n$ is topologically conjugate to an \emph{orthogonal}
  action for all $x\in\orbify{O}$.  That is, for each
  $x\in\orbify{O}$, there exists a (faithful) representation
  $\rho_x:\Gamma_x\to O(n)$, the orthogonal group, such that if
  $\gamma\cdot y$ denotes the $\Gamma_x$ action on $\tilde U_x$, there
  exists a homeomorphism $h$ of $\tilde U_x$ such that
  $h\circ(\gamma\cdot y)=[\rho_x(\gamma)](h(y))$ for all $y\in\tilde
  U_x$. By standard results, \cite[lemma~4.7.1]{MR88k:53002}, the
  class of locally smooth orbifold remains unchanged if we replace
  $O(n)$ by the general linear group, $GL(n)$, in our definition.
\end{definition}

\begin{definition}\label{SmoothOrbifold} Let $0 \le r \le \infty$.  An
  orbifold \orbify{O} is a \emph{smooth $C^r$ orbifold } if each
  $\Gamma_x$ acts by $C^r$ diffeomorphisms on $\tilde{U}_x$ and each
  embedding $\tilde\psi_{zx}$ is $C^r$. When $r=0$, a smooth $C^0$ orbifold
  is understood to be locally smooth.
\end{definition}

\begin{proposition}\label{SmoothIsLocallySmooth}
  If \orbify{O} is a smooth $C^r$ orbifold with $r>0$, then it is
  locally smooth. Moreover, the action of the local isotropy groups is
  smoothly $C^r$ conjugate to an orthogonal action.
\end{proposition}
\begin{proof} Let $\Gamma_x$ be the isotropy group of $x$, $U_x$ a
  neighborhood of $x$ with corresponding neighborhood $\tilde{U}_x$ of
  $0$ in $\R^n$ and homeomorphism $\phi_x:\tilde{U}_x/\Gamma_x\to U_x$
  with $\phi_x(0)=x$. By assumption, $\Gamma_x$ acts by $C^r$
  diffeomorphisms on $\tilde U_x$.  We denote the action of $\Gamma_x$
  by $(\gamma, \tilde{y}) \to \gamma\cdot \tilde{y}$ for all $\gamma
  \in \Gamma_x$ and $ \tilde{y} \in \tilde{U}_x$.  Note that
  $\Gamma_x\cdot 0=0$. Let $L_\gamma:T_0\tilde{U}_x \to
  T_0\tilde{U}_x$ be the linearization at $0$ of $\tilde{y} \to
  \gamma\cdot \tilde{y}$.  Note that $L_\gamma$, being the
  linearization {\em at} $0$, is a {\em fixed} linear map, and is
  therefore $C^\infty$. Define $F:\tilde {U}_x\to\R^n$ by
  $$F(\tilde{y}) = \frac{1}{|\Gamma_x|}\sum_{\eta \in \Gamma_x}L_\eta
  (\eta^{-1}\cdot \tilde{y})$$
  Then $F$ is $C^r$ since $L_\eta$ is $C^\infty$ and the action of
  $\Gamma_x$ is by $C^r$ diffeomorphisms. Also, $dF(0) = \hbox{Id}$
  and $F(\gamma \cdot \tilde{y}) = L_\gamma(F(\tilde{y}))$.  To see
  the last statement, note that
  \begin{align*}
    F(\gamma\cdot \tilde{y})& =\frac{1}{|\Gamma_x|}\sum_{\eta \in \Gamma_x}L_\eta (\eta^{-1}\gamma\cdot \tilde{y})\\
    &=\frac{1}{|\Gamma_x|}\sum_{\eta \in \Gamma_x}L_\eta ((\gamma^{-1}\eta)^{-1}\cdot \tilde{y})\\
    &=\frac{1}{|\Gamma_x|}\sum_{\mu\in \Gamma_x}L_{\gamma\mu} (\mu^{-1}\cdot \tilde{y}) && \text{where }\mu=\gamma^{-1}\eta\\
    &=\frac{1}{|\Gamma_x|}\sum_{\mu\in \Gamma_x}L_{\gamma}(L_\mu (\mu^{-1}\cdot \tilde{y}))\\
    &=L_\gamma\big(\frac{1}{|\Gamma_x|}\sum_{\mu\in \Gamma_x}L_\mu
    (\mu^{-1}\cdot \tilde{y})\big)=L_\gamma(F(\tilde{y}))
  \end{align*}

  So by the inverse function theorem, there is a neighborhood
  $\tilde{V}_x$ of $0$ in $\tilde{U}_x$ on which $F$ is an equivariant
  $C^r$ diffeomorphism. Thus, $F$ conjugates the action of $\Gamma_x$
  to the linear action $L_\gamma$ which in turn is linearly conjugate
  to an orthogonal action which we denote by $\rho_x(\gamma)$.
  $\rho_x$ is the required representation making $\orbify{O}$ locally
  smooth.
\end{proof}

\begin{definition}\label{orbifoldchart} 
  An \emph{orbifold chart} about $x$ in a (locally) smooth orbifold
  $\orbify{O}$ is a 4-tuple $(\tilde U_x, \Gamma_x, \rho_x, \phi_x)$
  where $\tilde U_x\cong\R^n$, $\Gamma_x$ is a finite group, $\rho_x$
  is a (faithful) representation of $\Gamma_x:\rho_x\in \Hom(\Gamma_x,
  O(n))$, and $\phi_x$ is a homeomorphism: $\phi_x:\tilde
  U_x/\rho_x(\Gamma_x)\to U_x$, where $U_x\subset X_{\orbify{O}}$ is a
  (sufficiently small) open relatively compact neighborhood of $x$,
  and $\phi_x(0)=x$.
\end{definition}

For convenience we will often refer to the neighborhood $U_x$ or
$(\tilde U_x, \Gamma_x)$ as an orbifold chart, and ignore the
representation $\rho_x$ and write $U_x=\tilde U_x/\Gamma_x$. If
necessary, we can assume that $\tilde U_x$ is an open metric ball in
$\R^n$ centered at the origin and denote by
$\pi_x:\tilde{U}_x\to\tilde{U}_x/\rho_x(\Gamma_x)$, the quotient map
defined by the action of $\rho_x(\Gamma_x)$ on $\tilde{U}_x$.

\begin{proposition}\label{FixedPointsAreSubmanifolds}
  Let $r\ge 0$. If \orbify{O} is a smooth $C^r$ orbifold then in each
  orbifold chart ${\tilde U}_x$ the fixed point set
  $\tilde{S}_x=\{\tilde{y}\in{\tilde U}_x\mid\Gamma_x\cdot
  \tilde{y}=\tilde{y}\}$ is a connected $C^r$ submanifold of ${\tilde
    U}_x$.
\end{proposition}
\begin{proof} Let $(\tilde U_x, \Gamma_x, \rho_x, \phi_x)$ be an
  orbifold chart about $x$. Since $\orbify{O}$ is $C^r$ smooth, the
  proof of proposition~\ref{SmoothIsLocallySmooth} gives the existence
  of $\Gamma_x$-equivariant $C^r$ diffeomorphism
  $F:\tilde{U}_x\to\R^n_{\rho_x}$, where $\R^n_{\rho_x}$ denotes
  $\R^n$ with the orthogonal $\Gamma_x$-action induced by the
  representation $\rho_x$. Thus, we have $F(\gamma\cdot
  \tilde{y})=[\rho_x(\gamma)](F(\tilde{y}))$. If $\tilde{y}\in
  \tilde{S}_x$, and $\tilde{z}=F(\tilde{y})$ then we have that
  $\tilde{z}=[\rho_x(\gamma)](\tilde{z})$, hence
  $F(\tilde{S}_x)\subset\bigcap_{\gamma\in\Gamma_x}\ker(\rho_x(\gamma)-I)$.
  Let $\tilde{W}=\bigcap_{\gamma\in\Gamma_x}\ker(\rho_x(\gamma)-I)$
  and let $\tilde{w}\in \tilde{W}$, with $F(\tilde{v})=\tilde{w}$ for
  some $\tilde{v}\in\tilde{U}_x$. Then
$$\tilde{v}=F^{-1}(\tilde{w})=F^{-1}[\rho_x(\gamma)](\tilde{w})=F^{-1}[\rho_x(\gamma)] F(\tilde{v})=%
F^{-1}F(\gamma\cdot \tilde{v})=\gamma\cdot \tilde{v}$$ for all
$\gamma\in\Gamma_x$.  Hence $\tilde{v}\in \tilde{S}_x$. We have shown
$F(\tilde{S}_x)=\tilde{W}$. Since $\tilde{W}$ is a subspace, we have
that $\tilde{S}_x=F^{-1}(\tilde{W})$ is a connected $C^r$ submanifold
of $\tilde{U}_x$.
\end{proof}

\subsection*{Stratification of an
  Orbifold}\label{StratificationSection}

\begin{definition}\label{Strata}
  Let $\orbify{O}$ be a connected $n$-dimensional locally smooth
  orbifold.  Given a point $x \in \orbify{O}$, there is a neighborhood
  $U_x$ of $x$ which is homeomorphic to a quotient $\tilde
  U_x/\Gamma_x$ where $\tilde U_x$ is homeomorphic to $\R^n$ and
  $\Gamma_x$ is a finite group acting orthogonally on $\R^n$. The
  definition of orbifold implies that the germ of this action in a
  neighborhood of the origin of $\R^n$ is unique.  We define the
  \emph{isotropy group of $x$} to be the group $\Gamma_x$.  The
  \emph{singular set}, $\Sigma_1$, of \orbify{O} is the set of points
  $x \in \orbify{O}$ with $\Gamma_x \ne \{ e\}$.
\end{definition}

We wish to define the notion of a {\em stratum} $\orbify{S}$ of
$\orbify{O}$. Roughly speaking, a stratum of $\orbify{O}$ is a maximal
connected subset $\orbify{S}$ of $\orbify{O}$ for which the $\Gamma_x$
action is constant for $x\in\orbify{S}$.  The formal definition is:

\begin{definition}\label{Stratum} Two points $x$, $y$ belong to the
  same stratum $\orbify{S}\subset\orbify{O}$ if there exists a chain
  of orbifold charts $\{U_x=U_0, U_1,\ldots, U_m=U_y\}$ so that for
  $0\le i\le m-1$ we have
  \begin{enumerate}
  \item $U_i\cap U_{i+1}\ne\emptyset$
  \item $\textup{Im}(\rho_i)=\textup{Im}(\rho_{i+1})$, and
  \item $\Gamma_i$ acts on $\tilde U_i\cap\tilde U_{i+1}$; that is,
    $\tilde U_i\cap\tilde U_{i+1}$ is $\Gamma_i$ invariant
  \end{enumerate}
  Here, $\rho_i\in\Hom(\Gamma_i,O(n))$ is the faithful representation
  of $\Gamma_i$ corresponding to the chart $U_i$.  By construction,
  the diagram below commutes (horizontal maps are simply inclusions):

\begin{equation*}
  \xymatrix{{\tilde U_i\cap\tilde U_{i+1}}\ar[rr]^{\subset}\ar[d]&&{\tilde U_{i+1}}\ar[d] && \\
    {(\tilde U_i\cap\tilde U_{i+1})/\Gamma_i}\ar[rr]^{\subset}\ar[d]%
    &&{\tilde U_{i+1}/\Gamma_{i+1}}\ar[d]\\
    {U_i\cap U_{i+1}}\ar[rr]^{\subset}&&{U_{i+1}} }
\end{equation*}
\end{definition}

It is easy to see that belonging to the same stratum is an equivalence
relation on $\orbify{O}$. Also, there can only be a finite number of
distinct strata on a compact orbifold. We have the following structure
result for strata:

\begin{proposition}\label{StratumNeighborhood} Let $\orbify{S}$ be a
  stratum of a smooth $C^r$ orbifold $\orbify{O}$. Then $\orbify{S}$
  is connected and there exists a connected smooth $C^r$ manifold
  $\tilde{U}$ and a $C^r$ action by a finite group $\Gamma$ on
  $\tilde{U}$ such that $\tilde{U}/\Gamma$ is a neighborhood of
  $\orbify{S}$ in $\orbify{O}$.
\end{proposition}
\begin{proof} From the definition of smooth orbifold we see that
  $\tilde{U}=\bigcup_{i=0}^m\tilde{U}_i$ inherits the structure of a
  connected smooth $C^r$ manifold. Let $\Gamma=\Gamma_0$ and
  $\rho=\rho_0$. By construction, we have an orthogonal action given
  by $\rho_0(\Gamma)$ of $\Gamma$ on $\tilde{U}$ and it is clear that
  $\tilde{U}/\Gamma$ is a neighborhood of $\orbify{S}$ in
  $\orbify{O}$. That $\orbify{S}$ is connected follows from
  proposition~\ref{FixedPointsAreSubmanifolds} and the fact that
  $\orbify{S}$ is the (continuous) projection of the fixed point
  subset
  $\tilde{S}=\{\tilde{u}\in\tilde{U}\mid\Gamma\cdot\tilde{u}=\tilde{u}\}$.
\end{proof}

\begin{definition}\label{StrataNeighborhoodChart} Let $\orbify{O}$ be
  a smooth $C^r$ orbifold. For $x\in\orbify{O}$, the stratum
  containing $x$ will be denoted by $\orbify{S}_x$. It is a
  suborbifold of $\orbify{O}$ (see definition~\ref{SubOrbifold}).  The
  corresponding $C^r$ manifold covering and finite group given in
  proposition~\ref{StratumNeighborhood} will be denoted by
  $\tilde{U}_{\orbify{S}_x}$ and $\Gamma_{\orbify{S}_x}$,
  respectively.  The neighborhood
  $\tilde{U}_{\orbify{S}_x}/\Gamma_{\orbify{S}_x}$ of $\orbify{S}_x$
  will be denoted by $U_{\orbify{S}_x}$ and the inverse image of
  $\orbify{S}_x$ in $\tilde{U}_{\orbify{S}_x}$ will be denoted by
  $\tilde{S}_x$.
\end{definition}

\subsection*{Products of Orbifolds}\label{OrbifoldProductSection}

Cartesian products of (locally) smooth orbifolds inherit a natural
(locally) smooth orbifold structure:

\begin{definition}\label{OrbifoldProduct}
  Let $\orbify{O}_i$ for $i = 1, 2$ be orbifolds.  The \emph{orbifold
    product} $\orbify{O}_1\times\orbify{O}_2$ is the orbifold having
  the following structure:
  \begin{enumerate}
  \item $X_{\orbify{O}_1\times\orbify{O}_2} = X_{\orbify{O}_1}\times
    X_{\orbify{O}_2}$.
  \item For each $(x_1, x_2)\in X_{\orbify{O}_1\times\orbify{O}_2}$
    and orbifold charts $U_i$ of $x_i$, $U_1\times U_2$ is an orbifold
    chart around $(x_1, x_2)$. Explicitly, $$(\tilde U_1\times\tilde
    U_2, \Gamma_{x_1}\times\Gamma_{x_2}, \rho_{x_1}\times\rho_{x_2},
    \phi_{x_1}\times\phi_{x_2})$$ is an orbifold chart around
    $(x_1,x_2)$.
  \end{enumerate}
  Note that the isotropy group $\Gamma_{(x_1, x_2)} =
  \Gamma_{x_1}\times\Gamma_{x_2}$.
\end{definition}

\subsection*{Suborbifolds}\label{SuborbifoldSection}

The definition of a suborbifold is somewhat more delicate than the
corresponding notion for a manifold.

\begin{definition}\label{SubOrbifold}
  A \emph{suborbifold} \orbify{P} of an orbifold \orbify{O} consists
  of the following.
  \begin{enumerate}
  \item A subspace $X_{\orbify{P}}\subset X_{\orbify{O}}$ equipped
    with the subspace topology
  \item For each $x\in X_{\orbify{P}}$ and neighborhood $W$ of $x$ in
    $X_{\orbify{P}}$ there is an orbifold chart $(\tilde U_x,
    \Gamma_x, \rho_x, \phi_x)$ about $x$ in \orbify{O} with
    $U_x\subset W$, a subgroup $\Lambda_x \subset \Gamma_x$ of the
    isotropy group of $x$ in \orbify{O} and a $\rho_x(\Lambda_x)$
    invariant vector subspace $\tilde V_x\subset \tilde U_x = \R^n$,
    so that $(\tilde V_x, \Lambda_x, \rho_x|_{\Lambda_x},\psi_x)$ is
    an orbifold chart for \orbify{P} and
  \item
    \begin{eqnarray*}
      V_x &=& \psi_x(\tilde V_x/\rho_x(\Lambda_x))\\
      &=& U_x\cap X_{\orbify{P}} \\
      &=& \phi_x(\pi_x(\tilde V_x))
    \end{eqnarray*}
    is an orbifold chart for $x$ in \orbify{P} where $\pi_x:\tilde
    U_x\to \tilde U_x/\rho_x(\Gamma_x)$ is the quotient map.
  \end{enumerate}

\end{definition}

\begin{remark}
  It is tempting to define the notion of an $m$--suborbifold
  $\orbify{P}$ of an $n$--orbifold $\orbify{O}$ simply by requiring
  $\orbify{P}$ to be locally modeled on $\R^m\subset\R^n$ modulo
  finite groups. That is, the local action on $\R^m$ is induced by the
  local action on $\R^n$. This is the definition adopted in
  \cite{Thurston78}.  It is equivalent to the added condition in our
  definition that $\Lambda_x = \Gamma_x$ at all $x$ in the underlying
  topological space of $\orbify{P}$. This more restrictive definition
  is not adequate for our needs as the following example shows.
\end{remark}

\begin{example}\label{Diagonal}
  Let \orbify{O} be a smooth $C^r$ orbifold. Let
  $\textup{diag}(\orbify{O})=\{(x,x)\mid x\in\orbify{O}\}\subset
  \orbify{O}\times\orbify{O}$ be the \emph{diagonal}.  Then
  $\textup{diag}(\orbify{O})$ is a suborbifold of
  $\orbify{O}\times\orbify{O}$ with isotropy group
  $\Gamma_{(x,x)}\cong\Gamma_x$ via the diagonal action $\gamma\cdot
  (\tilde x,\tilde x)=(\gamma\cdot\tilde x,\gamma\cdot\tilde x)$. See
  proposition~\ref{OrbifoldGraph}.  If we had chosen the more
  restrictive definition of suborbifold given in the last remark, then
  $\textup{diag}(\orbify{O})$ would {\em not} have been a suborbifold.
  For example, consider the orbifold $\R/\Z_2$ where $\Z_2$ acts on
  $\R$ via $\gamma\cdot x=-x$. The underlying topological space
  $X_{\orbify{O}}$ of $\orbify{O}$ is $[0, \infty)$ and the isotropy
  subgroups are $\{1\}$ for $x\in (0, \infty)$ and $\Z_2$ for $x = 0$.
  The isotropy subgroup of $(0,0)\in\R/\Z_2\times\R/\Z_2$ is
  $\Z_2\times\Z_2$, whereas the isotropy subgroup of $(0,0)$ in the
  diagonal suborbifold
  $\textup{diag}(\R/\Z_2)\subset\R/\Z_2\times\R/\Z_2$ must be
  isomorphic to $\Z_2$, as $\textup{diag}(\R/\Z_2)$ is a 1-dimensional
  suborbifold.
\end{example}

\begin{remark} 
  Let $\orbify{P}\subset\orbify{O}$ be a suborbifold.  Note that even
  though a point $p\in X_{\orbify{P}}$ may be in the singular set of
  \orbify{O}, it need not be in the singular set of \orbify{P}.
\end{remark}

\section{Orbifold Maps}\label{OrbifoldMapSection}
Intuitively, an orbifold map should be a map between
underlying topological spaces that has local lifts, but unfortunately
axiomatizing such a simple idea has proven difficult if one wants to provide a
definition that is very flexible.
We now discuss one such natural definition of maps between orbifolds. This
definition will elaborate on the definition that was given in the
paper \cite{MR2003g:58013}. In that paper, these maps were referred to
as unreduced orbifold maps because we distinguished among different
liftings of the same map of underlying topological spaces. From now on, we will refer to such maps
simply as orbifold maps.
In \cite{MR2003146}, our definition of (reduced) orbifold map did not distinguish
among different liftings. We will retain the term {\em reduced} for orbifold maps 
for which the particular choice of local lifts is ignored. Thus, a reduced orbifold map agrees with
the notion of orbifold map given in \cite[Def.~1.3]{AdemLeidaRuan}.

In what follows we use the notation given in
definitions~\ref{orbifold}, \ref{orbifoldchart} and
\ref{StrataNeighborhoodChart}.

\begin{definition}\label{OrbiMap}  
  A $C^0$ \emph{orbifold map} $(f,\{\tilde f_x\})$ between locally
  smooth orbifolds $\orbify{O}_1$ and $\orbify{O}_2$ consists of the
  following:
  \begin{enumerate}
  \item A continuous map $f:X_{\orbify{O}_1}\to X_{\orbify{O}_2}$ of
    the underlying topological spaces.
  \item For each $y\in {\mathcal{S}_x}$, a group homomorphism
    $\Theta_{f,y}:\Gamma_{\mathcal{S}_x}\to\Gamma_{f(y)}$.
  \item A $\Theta_{f,y}$-equivariant lift $\tilde f_y:\tilde
    U_y\subset\tilde{U}_{\mathcal{S}_x}\to\tilde V_{f(y)}$ where
    $(\tilde U_y,\Gamma_{\mathcal{S}_x}, \rho_y, \phi_y)$ is an
    orbifold chart at $y$ and $(\tilde V_{f(y)},\Gamma_{f(y)},
    \rho_{f(y)}, \phi_{f(y)})$ is an orbifold chart at $f(y)$.  That
    is, the following diagram commutes:
    \begin{equation*}
      \xymatrix{{\tilde U_y}\ar[rr]^{\tilde f_y}\ar[d]&&{\tilde V_{f(y)}}\ar[d]\\
        {\tilde U_y}/\Gamma_{\mathcal{S}_x}\ar[rr]^>>>>>>>>>>%
        {{\tilde f_y}/\Theta_{f,y}(\Gamma_{\mathcal{S}_x})}\ar[dd]&&{\tilde
          V_{f(y)}}/\Theta_{f,y}(\Gamma_{\mathcal{S}_x})\ar[d]\\
        &&{\tilde V_{f(y)}}/\Gamma_{f(y)}\ar[d]\\
        U_y\subset U_{\mathcal{S}_x}\ar[rr]^{f}&&V_{f(y)}
      }
    \end{equation*}
  \item (Equivalence) Two orbifold maps $(f,\{\tilde f_x\})$ and
    $(g,\{\tilde g_x\})$ are considered equivalent if for each
    $x\in\orbify{O}_1$, $\tilde f_x=\tilde g_x$ as germs. That is,
    there exists an orbifold chart $(\tilde U_x,\Gamma_x)$ at $x$ such
    that ${\tilde f_x}\vert_{\tilde U_x}={\tilde g_x}\vert_{\tilde
      U_x}$. Note that this implies that $f=g$.
  \end{enumerate}
\end{definition}

\begin{remark} Note that equivalence of two orbifold maps does not
  require that $\Theta_{f,x}=\Theta_{g,x}$. To see that this is
  justifiable, consider the example where $\orbify{O}$ is the orbifold
  $\R/\Z_2$ where $\Z_2$ acts on $\R$ via $x \to -x$ and $f$ is the
  constant map $f\equiv 0$. The underlying topological space
  $X_{\orbify{O}}$ of $\orbify{O}$ is $[0, \infty)$ and the isotropy
  subgoups are trivial for $x\in (0, \infty)$ and $\Z_2$ for $x = 0$.
  The map $\tilde f_0\equiv 0$ is a local equivariant lift of $f$ at
  $x=0$ using either of the homomorphisms $\Theta_{f,0}=\textup{Id}$
  or $\Theta'_{f,0}=\{e\}$. We do not wish to consider these as
  distinct orbifold maps.
\end{remark}

For convenience, we will often denote an orbifold map $(f,\{\tilde
f_x\})$ simply by $f$.

\begin{definition}\label{OrbiMapSmooth} An orbifold map
  $f:\orbify{O}_1\to\orbify{O}_2$ of $C^r$ smooth orbifolds is {\em
    $C^r$ smooth} if each of the local lifts $\tilde f_x$ may be
  chosen to be $C^r$.
\end{definition}

The next lemma is a technical result that states that a local lift $\tilde f_x$ chosen on a particular orbifold chart about $x$ uniquely specifies a local lift on any other orbifold chart about $x$. Hence, in definition~\ref{OrbiMap}, the $\tilde f_x$'s, once chosen, are independent of the choice of local charts.

\begin{lemma}\label{OrbiMapLiftExtenExistUniq}
  Let $f:\orbify{O}_1\to\orbify{O}_2$ be a $C^r$ orbifold map,
  $x\in\orbify{O}_1$, $U_x\subset W_x$ connected orbifold charts
  around $x$ and $V_{f(x)}\subset Z_{f(x)}$ connected orbifold charts
  around $f(x)$ in $\orbify{O}_2$ with $f(U_x)\subset V_{f(x)}$ and
  $f(W_x)\subset Z_{f(x)}$.  If $\tilde f_{U_x}$ is a lift of $f$ to
  $\tilde U_x$, then there is a unique lift $\tilde f_{W_x}$ of $f$ to
  $\tilde W_x$ extending $\tilde f_{U_x}$.
\end{lemma}

\begin{proof}
  Let $\tilde D_x\subset \tilde W_x$ and $\tilde D_{f(x)}\subset
  \tilde Z_{f(x)}$ be Dirichlet fundamental domains for the actions of
  the isotropy groups $\Gamma_x$ and $\Gamma_{f(x)}$ on $\tilde W_x$
  and $\tilde Z_{f(x)}$ respectively.  Then, $\tilde D_x\cap\tilde
  U_x$ and $\tilde D_{f(x)}\cap\tilde V_{f(x)}$ are also Dirichlet
  fundamental domains for the actions of the respective isotropy
  groups on $\tilde U_x$ and $\tilde V_{f(x)}$ respectively. Let
  $\tilde y\in\tilde U_x\cap\tilde D_x$ be a point in the non-singular
  set of $\orbify{O}_1$.  Without loss of generality, we may take
  $\tilde D_{f(x)}$ to be the Dirichlet fundamental domain containing
  $\tilde f_{\tilde U_x}(\tilde y)$ and so for any $\tilde z\in \tilde
  D_x$, there is a unique $\tilde w\in\tilde D_{f(x)}$ with
  $\pi_{f(x)}(\tilde w)= f(\pi_x(\tilde z))$. Now define the extension
  $\tilde f_{\tilde W_x}:\tilde W_x\to \tilde Z_{f(x)}$ via:
  \begin{equation*}
    \tilde f_{\tilde W_x}(\gamma\cdot\tilde z) = \Theta_{f,x}(\gamma)\cdot \tilde w
  \end{equation*}
  Uniqueness and continuity of the extension follow from the
  properties of Dirichlet domains.
\end{proof}

Given two orbifolds $\orbify{O}_i$, $i = 1,2$, the class of $C^r$
orbifold maps from $\orbify{O}_1$ to $\orbify{O}_2$ will be denoted by
$C^r_{\Orb}(\orbify{O}_1, \orbify{O}_2)$. If
$\orbify{O}_1=\orbify{O}_2=\orbify{O}$, we use the notation
$C^r_{\Orb}(\orbify{O})$ instead. The following was stated as a
proposition without proof in \cite{MR2003g:58013}.

\begin{example}[Lifts of the Identity Map]\label{IdentityMap}  Consider the identity map $\textup{Id}:\orbify{O}\to\orbify{O}$. Let 
  $x\in\orbify{O}$ and $(\tilde U_x,\Gamma_x)$ be an orbifold chart at
  $x$. From the definition of orbifold map, it follows (since
  $\Gamma_x$ is finite) that there exists $\gamma\in\Gamma_x$ such
  that a lift $\widetilde{\textup{Id}}_x:\tilde U_x\to\tilde U_x$ is
  given by $\widetilde{\textup{Id}}_x(\tilde y)=\gamma\cdot\tilde y$
  for all $\tilde y\in\tilde U_x$. Since $\widetilde{\textup{Id}}_x$
  is $\Theta_{\textup{Id},x}$ equivariant we have for
  $\delta\in\Gamma_x$:
  \begin{alignat*}{2}\widetilde{\textup{Id}}_x(\delta\cdot\tilde y) &
    = %
\Theta_{\textup{Id},x}(\delta)\cdot\widetilde{\textup{Id}}_x(\tilde y) &\quad & \text{hence }\\
\gamma\delta\cdot\tilde y & = \Theta_{\textup{Id},x}(\delta)\gamma\cdot\tilde y & & \text{which implies}\\
& & & \text{since $\Gamma_x$ acts effectively that}\\
\gamma\delta & = \Theta_{\textup{Id},x}(\delta)\gamma & & \text{or, equivalently,}\\
\Theta_{\textup{Id},x}(\delta) & = \gamma\delta\gamma^{-1}
\end{alignat*}
Thus, $\Theta_{\textup{Id},x}$ is an isomorphism of $\Gamma_x$, in
fact, an inner automorphism. Since two inner automorphisms,
$I_{\gamma_i}(\delta)=\gamma_i\delta\gamma_i^{-1}$, give rise to the same
automorphism of $\Gamma_x$ precisely when $\gamma_1=\zeta\gamma_2$
where $\zeta\in\textup{Center}(\Gamma_x)$, the number of possible
distinct choices for the homomorphism $\Theta_{\textup{Id},x}$ is
$\displaystyle\frac{|\Gamma_x|}{|\textup{Center}(\Gamma_x)|}$. In
particular, if $x$ is non--singular, or more generally, if $\Gamma_x$
is abelian, $\Theta_{\textup{Id},x}$ is the identity isomorphism on
$\Gamma_x$, and the identity map has exactly $|\Gamma_x|$ local lifts
over $x$.  Moreover, we see that the identity map between $C^r$
orbifolds is $C^r$. In fact, it is an example of a $C^r$ orbifold diffeomorphism (definition~\ref{homeos-and-diffeos}).
\end{example}

\begin{example}\label{TopMapsFromOrbifoldsToManifolds}
  Let $\orbify{O}$ be an orbifold and $X_\orbify{O}$ its underlying
  topological space. Let $N$ be a manifold or manifold with boundary
  (with trivial orbifold structure). Let $f:X_\orbify{O}\to N$ be a
  (topologically) continuous map; that is $f\in C^0(X_\orbify{O},N)$.
  Then $f$ is naturally an orbifold continuous map; that is $f\in
  C^0_{\Orb}(\orbify{O},N)$. To see this, note that since $N$ is a
  trivial orbifold, $\Gamma_{f(x)}=\{e\}$ for all $x\in\orbify{O}$.
  Thus, $\Theta_{f,x}$ is the constant homomorphism $\gamma\mapsto e$.
  Therefore, equivariant local lifts $\tilde f_x:\tilde U_x\to \tilde
  V_{f(x)}=V_{f(x)}$ may be defined via $\tilde f_x(\tilde
  y)=f\circ\pi_x(\tilde y)$ for $\tilde y\in \tilde U_x$. By
  construction $\tilde f$ is well-defined, continuous and unique, and
  thus $f\in C^0_{\Orb}(\orbify{O},N)$.
\end{example}

\begin{example}\label{MapsFromManifoldsToOrbifolds}
  Let $\orbify{O}$ be a smooth orbifold and let $N$ be a smooth
  manifold or manifold with boundary (with trivial orbifold
  structure). If $f\in C^r_{\Orb}(N,\orbify{O})$, then since
  $\Gamma_x=\{e\}$ for all $x\in N$ the homomorphism
  $\Theta_{f,x}:\Gamma_x\to\Gamma_{f(x)}$ is just $e\mapsto e$. Thus
  $f$ is merely a map from $N$ to $\orbify{O}$ with choice of local
  $C^r$ lifts. In the case where $\partial N\ne\emptyset$, this means
  that a local lift is $C^r$ over $N-\partial N$ with continuous
  extension to $\partial N$.
\end{example}

\begin{proposition}\label{OrbifoldGraph}
  Let $f\in C^r_{\Orb}(\orbify{O}_1,\orbify{O}_2)$, then the graph of
  $f$, $\mbox{graph}(f)$, defined by
  $$\textup{graph}(f) = \{(x, f(x))\in
  \orbify{O}_1\times\orbify{O}_2\}\subset\orbify{O}_1\times\orbify{O}_2
  $$
  is a $C^r$ suborbifold.  Note the isotropy group
  $\Gamma_{(x,y)}\cong\Gamma_x$ is acting on $\tilde U_x\times\tilde
  V_y$, a chart in $\orbify{O}_1\times\orbify{O}_2$, via the twisted
  diagonal action $\gamma\cdot (\tilde x,\tilde y)=(\gamma\cdot\tilde
  x,\Theta_{f,x}(\gamma)\cdot\tilde y)$.
\end{proposition}

\begin{proof} Let $x\in\orbify{O}_1$, $(\tilde U_x, \Gamma_x)$ a chart
  at $x$, $\Theta_{f,x}\in\textup{Hom}(\Gamma_x,\Gamma_{f(x)})$,
  $(\tilde V_{f(x)}, \Gamma_{f(x)})$ a chart at $f(x)$ and equivariant
  lift $\tilde f_x:\tilde U_x\to\tilde V_{f(x)}$ of $f$. That is,
  $\Theta_{f,x}(\gamma)\cdot\tilde f(\tilde x')=\tilde
  f(\gamma\cdot\tilde x')$ for all $\gamma\in\Gamma_x$ and $\tilde
  x'\in\tilde U_x$. For
  $(x,f(x))\in\textup{graph}(f)\subset\orbify{O}_1\times\orbify{O}_2$
  we have $\Gamma_{(x,f(x))}=\Gamma_x\times\Gamma_{f(x)}$. We need to
  give a suborbifold structure for $\textup{graph}(f)$.  Define the
  subgroup

  $\Gamma_\Theta=\{(\gamma,\Theta_{f,x}(\gamma))\mid\gamma\in\Gamma_x\}\subset%
  \Gamma_x\times\Gamma_{f(x)}$ and let $\tilde W_x=\{(\tilde x',\tilde
  f(\tilde x'))\mid\tilde x'\in\tilde U_x\}\subset\tilde
  U_x\times\tilde V_{f(x)}$.  Note that $\tilde W_x$ is
  $\Gamma_\Theta$ invariant: Suppose $(\tilde x',\tilde f(\tilde
  x'))\in\tilde W_x$ and
  $\delta=(\gamma,\Theta_{f,x}(\gamma))\in\Gamma_\Theta$. Then
$$\delta\cdot\left(\tilde x',\tilde f(\tilde x')\right)=%
\left(\gamma\cdot\tilde x',\Theta_{f,x}(\gamma)\cdot\tilde f(\tilde
  x')\right)=%
\left(\gamma\cdot\tilde x',\tilde f(\gamma\cdot\tilde x')\right)\in\tilde W_x$$
Thus, $\left(\tilde U_x\times\tilde V_{f(x)},\Gamma_x\times\Gamma_{f(x)},\rho_x\times\rho_{f(x)},\phi_x\times\phi_{f(x)}\right)$ is an orbifold chart around $(x,f(x))$ with 
$\left(\tilde W_x,{\Gamma_\Theta,\rho_x\times\rho_{f(x)}}\big\rvert_{\Gamma_\Theta},%
\psi_x=\phi_x\times\phi_{f(x)}\big\rvert_{\textup{graph}(f)}\right)$ the required suborbifold chart around 
$(x,f(x))\in\textup{graph}(f)$.
\end{proof}

\begin{definition}\label{homeos-and-diffeos}
  For any topological space, let $\Homeo(X)$ denote its group of
  homeomorphisms. For a $C^0$ orbifold $\orbify{O}$, denote by
  $\Homeo_{\Orb}(\orbify{O})$ the subgroup of $\Homeo(X_{\orbify{O}})$
  with $f,f^{-1}\in C^0_{\Orb}(\orbify{O})$.  If $\orbify{O}$ is a
  $C^r$ orbifold, $\Orbdiff^r({\orbify{O}})$, the {\em $C^r$ orbifold
    diffeomorphism group}, is the subgroup of
  $\Homeo_{\Orb}(\orbify{O})$ with $f,f^{-1} \in
  C^r_{\Orb}(\orbify{O})$.
\end{definition}

\begin{example} Consider the case of a so-called $\mathbb{Z}_p$-football $\orbify{O}=S^2/\mathbb{Z}_p$ where $\mathbb{Z}_p$ acts on $S^2\subset\mathbb{R}^3$ by rotation about the $z$-axis by an angle $2\pi/p$. It is an example, in the language of Thurston, of a {\em good} orbifold $\orbify{O}=M/\Gamma$ where $M$ is a smooth manifold and $\Gamma$ acts effectively on $M$ as a proper discontinuous group of diffeomorphisms on $M$. This type of orbifold is referred to as an {\em effective global quotient} in \cite{AdemLeidaRuan}. There are two singular points corresponding to the north and south poles.
Let $\mathscr{ID}$ denote the subgroup of $\Orbdiff^r({\orbify{O}})$ comprised of all lifts of the identity map. 
Then $\mathscr{ID}\cong\mathbb{Z}_p\times\mathbb{Z}_p$. If we let $\Diff^r_{\mathbb{Z}_p}(M)\subset\Orbdiff^r({\orbify{O}})$ denote the (global)
$\mathbb{Z}_p$-equivariant diffeomorphisms of $M$ and let $\mathscr{ID}_{\mathbb{Z}_p}\subset\Diff^r_{\mathbb{Z}_p}(M)$ denote the 
$\mathbb{Z}_p$-equivariant lifts of the identity, then $\mathscr{ID}_{\mathbb{Z}_p}\cong\mathbb{Z}_p$. This example shows that, in general, 
$\Orbdiff^r({\orbify{O}})$ will be strictly larger than $\Diff^r_{\Gamma}(M)$ for a good orbifold $\orbify{O}=M/\Gamma$.
\end{example}

Recall the following terminology \cite{MR0448362}: Let $\mathscr{R}$
be a $C^r$ smooth structure on an orbifold $\orbify{O}$. A $C^s$
smooth structure $\mathscr{S}$ on $\orbify{O}$, $s>r$, is {\em
  compatible with} $\mathscr{R}$ if $\mathscr{S}\subset\mathscr{R}$.
This means that orbifold charts in $(\orbify{O},\mathscr{S})$ are
orbifold charts in $(\orbify{O},\mathscr{R})$ in the sense that the
identity map of $\orbify{O}$ is a element of
$\Orbdiff^r({\orbify{O}})$.  As in the classical case of smooth
manifolds \cite{MR1503303}, we have the following result on raising
the differentiablity of smooth orbifold structures.

\begin{proposition}\label{OrbifoldsAreSmoothable}
  Let $\mathscr{R}$ be a $C^r$ smooth structure on an orbifold
  $\orbify{O}$, $r\ge 1$. For every $s$, $r<s\le\infty$, there exists
  a compatible $C^s$ smooth structure $\mathscr{S}\subset\mathscr{R}$,
  and $\mathscr{S}$ is unique up to $C^s$ orbifold diffeomorphism.
\end{proposition}

\begin{proof} In light of definition~\ref{SmoothOrbifold} and
  example~\ref{IdentityMap}, one merely need use the results of
  Palais~\cite{MR0268912}.
\end{proof}

\section{Function Space Topologies}\label{CrTopologySection}

In this section, we assume that $\orbify{O}_i$ are smooth $C^r$
orbifolds and define the (strong/fine/Whitney) $C^r$ topology on
$C^r_{\Orb}(\orbify{O}_1, \orbify{O}_2)$.  For $f\in
C^r_{\Orb}(\orbify{O}_1, \orbify{O}_2)$, we first define a $C^0$
neighborhood of $f$ and corresponding $C^0$ topology on
$C^r_{\Orb}(\orbify{O}_1, \orbify{O}_2)$. Although we will introduce a
Riemannian structure later, for our purposes now we make the
observation that orbifolds are metrizable: Just let $U=\tilde
U/\Gamma=\pi(U)$ be any orbifold chart of $\orbify{O}$. Since $\Gamma$
is finite, we may define a metric on $U$ by $d_U(x,y)=d_{\tilde
  U}\left(\pi^{-1}(x),\pi^{-1}(y)\right)$ where $d_{\tilde U}$ is the
usual Euclidean metric on $\tilde U$. This makes $\orbify{O}$ locally
metrizable. Since all orbifolds are assumed paracompact and Hausdorff,
the Smirnov metrization theorem \cite{MR0464128} implies $\orbify{O}$
is metrizable and second countable.

\begin{definition}\label{C0-topology}
  Let $f:\orbify{O}_1\to\orbify{O}_2$ be a $C^r$ orbifold map. Let
  $\mathcal{C} = \{C_i\}$ be a locally finite covering of
  $\orbify{O}_1$ by relatively compact, open sets such that $\overline
  C_i\subset U_i$ and $f(\overline C_i)\subset V_i$ where $U_i$ and
  $V_i$ are (open) relatively compact orbifold charts.  Let
  $\{\varepsilon_i\}$ be a collection of positive constants. Let
  $\mathscr{N}^0(f,\varepsilon_i;\mathcal{C})$ consist of all $g\in
  C^r_{\Orb}(\orbify{O}_1, \orbify{O}_2)$ such that for all $i$,
  $g(C_i)\subset V_i$ and $\|\tilde f_x(\tilde y)-\tilde g_x(\tilde
  y)\|_{\tilde V_i}<\varepsilon_i$ for all $x\in C_i$ and $\tilde
  y\in\pi_x^{-1}(C_i\cap U_x)$.  The sets
  $\mathscr{N}^0(f,\varepsilon_i;\mathcal{C})$ form a neighborhood base for a
  topology on $C^r_{\Orb}(\orbify{O}_1, \orbify{O}_2)$, which we call
  the (orbifold) {\em $C^0$ topology} relative to $\mathcal{C}$ and we
  refer to $C^r_{\Orb}(\orbify{O}_1, \orbify{O}_2)$ with this topology
  as $C^r_{\Orb}(\orbify{O}_1, \orbify{O}_2;\mathcal{C})$.
\end{definition}

To define the (strong/fine/Whitney) $C^s$ topology on
$C^r_{\Orb}(\orbify{O}_1, \orbify{O}_2)$ for $1\le s\le r$, we simply
require, in addition, that local lifts are $C^s$ close in the usual
$C^s$ topology. In particular we have,

\begin{definition}\label{Cs-topology}
  Let $f:\orbify{O}_1\to\orbify{O}_2$ be a $C^r$ orbifold map. Define
  $\mathscr{N}^s(f,\varepsilon_i;\mathcal{C})$ to be those maps
  $g\in\mathscr{N}^0(f,\varepsilon_i;\mathcal{C})$ such that for all $1\le k\le
  s$, $\|\partial^k\tilde f_x(\tilde y)-\partial^k\tilde g_x(\tilde
  y)\|<\varepsilon_i$ for all $x\in C_i$ and $\tilde
  y\in\pi_x^{-1}(C_i\cap U_x)$.  This means that the local lifts of
  $f$ and $g$ have all partial derivatives of order $\le s$ within
  $\varepsilon_i$ at each point of $\tilde y\in\pi_x^{-1}(C_i\cap
  U_x)$.  Sets of this type form a neighborhood base for the
  (orbifold) $C^s$ {\em topology} on $C^r_{\Orb}(\orbify{O}_1,
  \orbify{O}_2)$ relative to the atlas $\mathcal{C}$.  The $C^\infty$
  {\em topology} relative to $\mathcal{C}$ on
  $C^\infty_{\Orb}(\orbify{O}_1, \orbify{O}_2)$ is defined to be the
  union of the topologies induced by the inclusion maps
  $C^\infty_{\Orb}(\orbify{O}_1,
  \orbify{O}_2;\mathcal{C})\hookrightarrow C^r_{\Orb}(\orbify{O}_1,
  \orbify{O}_2;\mathcal{C})$ for finite $r$ and as above, and
  $C^\infty_{\Orb}(\orbify{O}_1, \orbify{O}_2)$ with this topology
  will be denoted by $C^\infty_{\Orb}(\orbify{O}_1,
  \orbify{O}_2;\mathcal{C})$ as above.
\end{definition}

\begin{remark} If both $\orbify{O}_1$ and $\orbify{O}_2$ are compact,
  then the coverings $\{C_i\}$ are finite and $\varepsilon_i$ may be
  chosen to be a constant $\varepsilon$ for all $i$. The resulting
  topologies induced by the neighborhood base
  $\mathscr{N}^s(f,\varepsilon)$ on $C^r_{\Orb}(\orbify{O}_1,
  \orbify{O}_2)$ are equivalent to the topologies in
  definitions~\ref{C0-topology} and \ref{Cs-topology} given above.
\end{remark}

\begin{proposition}\label{IndependenceOfCover}
  The topology on $C^r_{\Orb}(\orbify{O}_1, \orbify{O}_2)$ is
  independent of the cover $\mathcal{C}$. That is, the spaces
  $C^r_{\Orb}(\orbify{O}_1, \orbify{O}_2;\mathcal{C})$ and
  $C^r_{\Orb}(\orbify{O}_1, \orbify{O}_2;\mathcal{C}')$ are
  homeomorphic for any two covers $\mathcal{C}$ and $\mathcal{C}'$ as
  in definition \ref{Cs-topology} and any value of $r$ where $0\le r\le\infty$.
\end{proposition}

The proof depends on the following lemma.  To aid both
the statement and proof of the following lemma, the following
notation will be useful.  For $f\in
C^r_{\Orb}(\orbify{O}_1,\orbify{O}_2)$, $U$ a chart about
$x\in\orbify{O}_1$, $V$ a chart about $f(x)\in\orbify{O}_2$ and relatively
compact connected open sets $x \in C'\subset\overline C'\subset C\subset\overline C\subset U$, define
\begin{align*}
  \mathscr{N}^s(f,\varepsilon;C) &= \{ g\in C^r_{\Orb}(\orbify{O}_1,\orbify{O}_2)\text{ such that } \\
  &\|\partial^k\tilde f(\tilde y)- \partial^k\tilde g(\tilde y)\| <
  \varepsilon\text{ for all } \tilde y\in\tilde C \text{ and all }
  k\le s \}
\end{align*}
\begin{align*}
  \mathscr{N}^s(f,\varepsilon;C,C')& = \{g\in \mathscr{N}^s(f,\varepsilon;C')\text{ such that }\\
  & \|\partial^k f(y) -\partial^k g(y)\|<\varepsilon\text{ for all } y\in C - \Sigma_1\text{ and }\\
  &\|f(y) - g(y)\|<\varepsilon\text{ for all }y\in C\}
\end{align*}

\begin{lemma}\label{EquivalenceUnderRefinement}
  Let $f$, $x$, $U$, $C$ and $C' \subset C$ be as above, then for each
  $\varepsilon > 0$ there is a $\delta > 0$ so that
  \begin{equation*}
    \mathscr{N}^r(f,\delta;C,C')\subset\mathscr{N}^r(f,\varepsilon;C)
  \end{equation*}
\end{lemma}

\begin{proof}
  The proof is by contradiction.  Assuming the contrary implies that
  there is an $\varepsilon> 0$ and a sequence $\{g_n\}\subset
  C^r_{\Orb}(\orbify{O}_1, \orbify{O}_2)$ so that
  \begin{equation*}
    g_n\in \mathscr{N}^r(f,2^{-n};C,C')\textup{ and } g_n\notin\mathscr{N}^r(f,\varepsilon;C)
  \end{equation*}
  For each $y\in C$, let $\Gamma_{f(y)}$ to be the isotropy group of
  $f(y)$ and $\theta_{f(y)f(x)}:\Gamma_{f(y)}\to\Gamma_{f(x)}$ the
  injective homomorphism of definition~\ref{orbifold}.  Let $N(x,y)$
  denote the index of $\theta_{f(y)f(x)}(\Gamma_{f(y)})$ in
  $\Gamma_{f(x)}$,
  $\left|\Gamma_{f(x)}:\theta_{f(y)f(x)}(\Gamma_{f(y)})\right|$ and
  let $\gamma_i$, $i = 1,\ldots, N(x, y)$ the corresponding coset
  representatives. Then there is a neighborhood, $\tilde V_{\tilde
    f(\tilde y)}$ of $\tilde f(\tilde y)$ in $\tilde V$ so that
  $\gamma_i\cdot\tilde V_{\tilde f(y)}\cap\,\gamma_j\cdot\tilde
  V_{\tilde f(\tilde y)} = \emptyset$ if $i\ne j$. Thus, the
  projection $\pi:\tilde V/\theta_{f(x)f(y)}(\Gamma_{f(y)})\to V$ is a
  local isometry over $\tilde V_{f(y)}$ by our choice of metric.  For
  any $\tilde y\in\tilde C$ let $\tilde W_{\tilde y} = \tilde
  f^{-1}\left(\tilde V_{\tilde f(\tilde y)}\right)$. $\left\{\tilde
    W_{\tilde y}\right\}$ is an open cover of $\tilde C$. Compactness
  of $\overline{\tilde C}$ yields a finite subcover $\tilde W_{\tilde
    y_1},\ldots,\tilde W_{\tilde y_M}$.  Without loss of generality,
  we may also uniformly bound the radii of the neighborhoods
  $V_{f(y)}$ in the range so that this cover is non-trivial.

  Now let $\tilde D\subset \tilde C$ be the maximal domain defined by
  \begin{equation*}
    \tilde D = \{\tilde z\in\tilde C\  |\  \tilde g_n(\tilde z)\to \tilde f(\tilde z)\text{ pointwise}\}
  \end{equation*}
  A Cantor diagonal argument shows that the limit point of any
  sequence $\tilde z_n\to\tilde z$ is also in $\tilde D$ and so
  $\tilde D$ is closed and therefore a compact set containing $\tilde
  C'$. Thus, there are points $\tilde y_{\alpha_1},\ldots,\tilde
  y_{\alpha_k}\subset\left\{\tilde y_1,\ldots,\tilde y_M\right\}$ so
  that $\tilde W_{\tilde y_{\alpha_1}},\ldots,\tilde W_{\tilde
    y_{\alpha_k}}$ cover $\tilde D$ and $\tilde D\cap\tilde W_{\tilde
    y_{\alpha_i}}\ne\emptyset$ for $i = 1,\ldots,k$.  By shrinking the
  $\tilde W_{\tilde y_{\alpha_i}}$'s we may assume that they still
  cover $\tilde D$ and they also satisfy $\tilde g_n(\tilde W_{\tilde
    y_{\alpha_i}}) \subset\tilde V_{\tilde f(\tilde y_{\alpha_i})}$
  for $n$ sufficiently large and all $i$.  Picking $\tilde z_i
  \in\tilde D\cap\tilde W_{\tilde y_{\alpha_i}}$ for each $i$ we have
  by definition of the $\tilde W$'s that
  \begin{equation*}
    \|\tilde g_n(\tilde z) - \gamma_i\cdot\tilde f(\tilde z)\| = \|g_n(z) - f(z)\|
  \end{equation*}
  for all $\tilde z\in\tilde W_{\tilde y_{\alpha_i}}$ and some coset
  representative $\gamma_i$ of $\theta_{f(x)f(y)}(\Gamma_{f(y)})$ in
  $\Gamma_{f(x)}$.  By evaluating at some $\tilde z_i\in\tilde
  D\cap\tilde W_{\tilde y_{\alpha_i}}$, the definition of $\tilde D$
  implies we must have $\gamma_i = e$ and thus, $\tilde g_n(\tilde
  z)\to\tilde f(\tilde z)$ for all $\tilde z\in\tilde W_{\tilde
    y_{\alpha_i}}$.  Since this holds for each $i=1,\ldots,k$, $\tilde
  g_n(\tilde z)\to\tilde f(\tilde z)$ for all $\tilde
  z\in\bigcup_{i=1}^k\tilde W_{\tilde y_{\alpha_i}}$ of which $\tilde
  D$ is a proper subset. This contradicts the maximality of $\tilde
  D$.
\end{proof}

  \begin{proof}[Proof of proposition~\ref{IndependenceOfCover}]
    Given two open covers $\mathcal{C}$ and $\mathcal{C}'$, take an
    open cover $\mathcal{C}''$ that refines them both.  Clearly the
    inclusion maps
  $$ C^r_{\Orb}(\orbify{O}_1,\orbify{O}_2;\mathcal{C})  \hookrightarrow C^r_{\Orb}(\orbify{O}_1,\orbify{O}_2;\mathcal{C}'')\text{ and }
  C^r_{\Orb}(\orbify{O}_1,\orbify{O}_2;\mathcal{C}') \hookrightarrow
  C^r_{\Orb}(\orbify{O}_1,\orbify{O}_2;\mathcal{C}'') $$ induced by
  restriction to the common refinement $\mathcal{C}''$ in each of the
  covers $\mathcal{C}$ and $\mathcal{C}'$ show that the topology on
  $C^r_{\Orb}(\orbify{O}_1,\orbify{O}_2;\mathcal{C}'')$ is coarser
  than either of the topologies induced by $\mathcal{C}$ or
  $\mathcal{C}'$. We now show that
  $C^r_{\Orb}(\orbify{O}_1,\orbify{O}_2;\mathcal{C}'')$ is, in fact,
  homeomorphic to $C^r_{\Orb}(\orbify{O}_1,\orbify{O}_2;\mathcal{C})$.
  
  Since sets of the form $\mathscr{N}^r(f,\varepsilon;C)$ for
  $C\in\mathcal{C}$ form a subbase for the topology of
  $C^r_{\Orb}(\orbify{O}_1,\orbify{O}_2;\mathcal{C})$, it suffices to
  find a neighborhood of $f$ in
  $C^r_{\Orb}(\orbify{O}_1,\orbify{O}_2;\mathcal{C}'')$ contained in
  $\mathscr{N}^r(f,\varepsilon;C)$.  Let
  $C_1'',\ldots,C_k''\in\mathcal{C}''$ be a cover of
  $C\in\mathcal{C}$.  For any $\delta > 0$
  \begin{equation*}
    \bigcap_{i=1}^k\mathscr{N}^r(f,\delta;C_k'')\subset\mathscr{N}^r(f,\delta;C,C_i'')
  \end{equation*}
  Therefore, by lemma~\ref{EquivalenceUnderRefinement},
  $\mathscr{N}^r(f,\varepsilon;C)$ is open in
  $C^r_{\Orb}(\orbify{O}_1,\orbify{O}_2;\mathcal{C}'')$ and thus we
  may conclude that
  $C^r_{\Orb}(\orbify{O}_1,\orbify{O}_2;\mathcal{C}'')$ and
  $C^r_{\Orb}(\orbify{O}_1,\orbify{O}_2;\mathcal{C})$ are
  homeomorphic.  Similarly,
  $C^r_{\Orb}(\orbify{O}_1,\orbify{O}_2;\mathcal{C}'')$ and
  $C^r_{\Orb}(\orbify{O}_1,\orbify{O}_2;\mathcal{C}')$ are
  homeomorphic. Thus,
  $C^r_{\Orb}(\orbify{O}_1,\orbify{O}_2;\mathcal{C})$ and
  $C^r_{\Orb}(\orbify{O}_1,\orbify{O}_2;\mathcal{C}')$ are
  homeomorphic as claimed.
\end{proof}

From now on, we drop the dependence of topology on
$C^r_{\Orb}(\mathcal{O}_1,\mathcal{O}_2;\mathcal{C})$ on the cover
$\mathcal{C}$, and will simply use the notation
$C^r_{\Orb}(\mathcal{O}_1,\mathcal{O}_2)$ for the set of orbifold
functions with the $C^r$ topology as in definition~\ref{Cs-topology}.
For the remainder of the paper, whenever function spaces between
orbifolds are mentioned, we will assume that the source orbifolds are
compact.

\begin{definition}\label{ConvergenceOfFunctions}
  For a fixed cover $\mathcal{C}$ by orbifold charts and any
  $\varepsilon > 0$, put
  \begin{equation*}
    \mathscr{N}^s(f,\varepsilon) = \{g\in C^r_{\Orb}(\orbify{O}_1,\orbify{O}_2)\  |\  g\in\mathscr{N}^s(f,\varepsilon;C)\mbox{ for all }C\in\mathcal{C}\}
  \end{equation*}
  As in the case for compact manifolds, for a compact orbifold
  $\orbify{O}_1$, we define for $f$ and $g\in
  C^r_{\Orb}(\orbify{O}_1,\orbify{O}_2)$ a distance
  \begin{equation*}
    d_s(f,g) = \inf\{\varepsilon > 0\  |\  f\in\mathscr{N}^s(g,\varepsilon)\mbox{ and }g\in\mathscr{N}^s(f,\varepsilon)\}
  \end{equation*}
  where the dependence on the orbifold atlas used has been supressed.
\end{definition}

\begin{remark}
  Compactness of $\orbify{O}_1$ implies (as in the usual manifold case) that the metric topology
  induced by the metric $d_s$ as above is equivalent to the $C^s$
  topology on $C^r_{\Orb}(\orbify{O}_1,\orbify{O}_2)$ given by the
  orbifold atlas $\mathcal{C}$ (and hence to the topology induced by
  any other atlas by proposition \ref{IndependenceOfCover}).
\end{remark}

\begin{proposition}\label{CrTopologyComplete}
  Let $\orbify{O}_i$ be compact $C^r$ orbifolds, $1\le r\le\infty$.
  For $1\le s \le r$, $C^r_{\Orb}(\orbify{O}_1, \orbify{O}_2)$ with
  the $C^s$ topology relative to $\mathcal{C}$ is a separable metric
  space. If $s=r$, then this metric space is complete.
\end{proposition}
\begin{proof}

  Let $\{f_n\}\subset C_{\Orb}^r(\orbify{O}_1,\orbify{O}_2)$ be a
  Cauchy sequence in the $C^r$ topology.  For any $x\in\orbify{O}_1$,
  orbifold charts $U_x$ about $x$ and $V\subset \orbify{O}_2$
  containing $\bigcup_n f_n(U_x)$, the lifts $\{\tilde f_n:\tilde
  U_x\to \tilde V\}$ are a sequence of $\Gamma_x$-equivariant
  functions converging uniformly in the $C^r$ topology on compact
  subsets of $\tilde U_x$.  Therefore they converge to a $C^r$,
  $\Gamma_x$-equivariant function $\tilde f:U_x\to\tilde V$ which is a
  lift of the function $f(x) = \lim f_n(x)$.  Thus, the limit function
  $f\in C^r_{\Orb}(\orbify{O}_1,\orbify{O}_2)$ which proves
  completeness.  For separability, note that for any $f\in
  C^r_{\Orb}(\orbify{O}_1,\orbify{O}_2)$, each lift $\tilde f_x:\tilde
  U_x\to\tilde V_{f(x)}$ may be approximated by a polynomial $\tilde
  g_x:\tilde U_x\to\tilde V_{f(x)}$.  To get a $\Gamma_x$-equivariant
  approximation by a polynomial we average $\tilde g_x$ over
  $\Gamma_x$. That is, we define $\tilde G_x:\tilde U_x\to\tilde
  V_{f(x)}$ by
  $$\tilde G_x(\tilde z)=\frac{1}{|\Gamma_x|}\sum_{\gamma\in\Gamma_x}\Theta_{f,x}(\gamma)\cdot\tilde g_x(\gamma^{-1}\cdot\tilde z)$$
  Since
  \begin{align*}
    \tilde G_x(\delta\cdot\tilde z)& =\frac{1}{|\Gamma_x|}\sum_{\gamma
      \in \Gamma_x}\Theta_{f,x}(\gamma)%
    \cdot\tilde g_x(\gamma^{-1}\delta\cdot\tilde z)\\
    &=\frac{1}{|\Gamma_x|}\sum_{\gamma \in \Gamma_x}\Theta_{f,x}(\delta)\Theta_{f,x}(\delta^{-1}\gamma)%
    \cdot\tilde g_x({(\delta^{-1}\gamma)}^{-1}\cdot\tilde z)\\
    &=\Theta_{f,x}(\delta)\cdot\frac{1}{|\Gamma_x|}\sum_{\mu \in \Gamma_x}\Theta_{f,x}(\mu)%
    \cdot\tilde g_x(\mu^{-1}\cdot\tilde z) & \text{where }\mu=\delta^{-1}\gamma\\
    &=\Theta_{f,x}(\delta)\cdot\tilde G_x(\tilde z)
  \end{align*}
  we see that $\tilde G_x$ satisfies the same equivariance relation as
  $\tilde f_x$ and thus $\tilde G_x\in
  C^r_{\Orb}(\orbify{O}_1,\orbify{O}_2)$. Since averaging is distance
  nonincreasing, we have produced an approximation of $\tilde f_x$ by
  $\Gamma_x$-equivariant polynomials. Furthermore, because there can
  be only finitely many lifts of $f$ over any orbifold chart,
  compactness of $\orbify{O}_1$ implies that the space
  $C^r_{\Orb}(\orbify{O}_1,\orbify{O}_2)$ is separable as the
  equivariant polynomials form a countable dense set.

\end{proof}

\section{The Tangent Orbibundle and its
  Sections}\label{TangentObibundleSection}

We now define the tangent orbibundle of a smooth $C^{r+1}$ orbifold.
It is a special case of the more general notion of a linear orbibundle
given in \cite{MR2003g:58013}.

\begin{definition}\label{tangentorbibundle}
  Let $\orbify{O}$ be an $n$--dimensional $C^{r+1}$ orbifold.  The
  \emph{tangent orbibundle} of $\orbify{O}$,
  $p:T\orbify{O}\to\orbify{O}$, is the $C^{r}$ orbibundle defined as
  follows.  If $(\tilde U_x,\Gamma_x)$ is an orbifold chart around
  $x\in \orbify{O}$ then $p^{-1}(U_x)\cong(\tilde
  U_x\times\R^n)/\Gamma_x$ where $\Gamma_x$ acts on $\tilde
  U_x\times\R^n$ via $\gamma\cdot(\tilde y, \tilde v) =
  \left(\gamma\cdot\tilde y, d\gamma_{\tilde y}(\tilde v)\right)$. In
  keeping with tradition, we denote the fiber $p^{-1}(x)$ over $x\in
  U_x$ by $T_x\orbify{O}\cong\R^n/\Gamma_x$. Note that, in general, if
  $\Gamma_x$ is non-trivial then $T_x\orbify{O}$ will be a convex cone
  rather than a vector space. Locally we have the diagram:
  
  \begin{equation*}
    \xymatrix{
      {\tilde U_x\times\R^n}\ar[rrr]^{\Pi_x}\ar[d]_{\text{pr}_1} & & &{(\tilde U_x\times\R^n)/\Gamma_x}\ar[d]^{p}\\
      {\tilde U_x}\ar[rrr]^{\pi_x}& &  & {U_x}
    }
  \end{equation*}
  where $\text{pr}_1:\tilde U_x\times\R^n\to\tilde U_x$ denotes the
  projection onto the first factor $(\tilde y,\tilde v)\mapsto\tilde
  y$ (which is a specific choice of lift of $p$).

\end{definition}

\begin{definition}\label{OrbiSection}
  A $C^r$ \emph{orbisection} of the tangent orbibundle $T\orbify{O}$
  is a $C^r$ orbifold map $\sigma:\orbify{O}\to T\orbify{O}$ such that
  $p\circ \sigma = \text{Id}_\orbify{O}$ and for any $x\in\orbify{O}$
  and chart $U_x$ about $x$, we have $\text{pr}_1\circ\tilde\sigma_x =
  \text{Id}_{\tilde U_x}$. In particular, it follows that
  $\Theta_{\sigma,x}=\text{Id}:\Gamma_x\to\Gamma_x$ and thus
  orbisections have unique equivariant lifts over orbifold charts.
\end{definition}

We have the following structure result which was first stated in
\cite{MR2003g:58013}.

\begin{proposition}\label{OrbiSecsAreVecSpace}
  The set $\mathscr{D}^r_{\Orb}(\orbify{O})$ of $C^r$ orbisections of
  the tangent orbibundle $T\orbify{O}$ is naturally a real vector
  space with the vector space operations being defined pointwise.
\end{proposition}

\begin{proof} Let $\sigma\in\mathscr{D}^r_{\Orb}(\orbify{O})$. Locally
  we have the diagram:

\begin{equation*}
  \xymatrix{
    {\tilde U_x}\ar[rrr]^{\tilde\sigma_x}\ar[d]_{\pi_x} & & &{\tilde U_x\times\R^n}\ar[d]^{\Pi_x}\\
    U_x\ar[rrr]^{\sigma_x}\ar[drrr]^{\textup{Id}}& &  & {p^{-1}(U_x)=(\tilde U_x\times\R^n)/\Gamma_x}\ar[d]^p\\
    &  &  & U_x
  }
\end{equation*}
and we can write for $y\in U_x$, $\sigma(y)=(y,s(y))$ where $s(y)\in
T_y\orbify{O}\cong\R^n/\theta_y(\Gamma_y)$ ($\theta_y$ is the
injective homomorphism which appears in definition~\ref{orbifold}).
Let $\tilde\sigma_x$ be the lift of $\sigma$.  Then
$\tilde\sigma_x(\tilde y)=(\tilde y,\tilde s(\tilde y))$, where
$\tilde s:\tilde U_x\to\R^n$ is such that $\tilde s(\delta\cdot\tilde
y)=d\delta_{\tilde y}(\tilde s(\tilde y))$. In particular, since
$\tilde x$ is a fixed point of the $\Gamma_x$ action on $\tilde U_x$,
we have $\tilde s(\tilde x)=\tilde s(\delta\cdot\tilde x)=%
d\delta_{\tilde x}(\tilde s(\tilde x))$.  Thus $\tilde s(\tilde x)$ is
a fixed point of the (linear) action of $\Gamma_x$ on $\R^n$. Note
that the set of such fixed points forms a {\em vector subspace} of
$\R^n$. As a result we may define a real vector space structure on
$\mathscr{D}^r_{\Orb}(\orbify{O})$ as follows: For
$\sigma_i\in\mathscr{D}^r_{\Orb}(\orbify{O})$, let
$\tilde\sigma_{i,x}$ be local lifts at $x$ as above. Define
\begin{align*}
  (\sigma_1+\sigma_2)(y) & =
  \Pi_x\big((\tilde\sigma_{1,x}+\tilde\sigma_{2,x})(\tilde y)\big)=%
  \Pi_x\big((\tilde y,\tilde s_{1}({\tilde y})+\tilde s_{2}({\tilde y}))\big)=\sigma_1(y)+\sigma_2(y)\\
  (\lambda\sigma)(y) & = \Pi_x\big((\lambda\tilde\sigma_{x})(\tilde
  y)\big)=%
  \Pi_x\big((\tilde y,\lambda\tilde s({\tilde
    y}))\big)=\lambda(\sigma(y))
\end{align*}

\end{proof}

In light of the previous proposition, we make the following

\begin{definition} Let $\orbify{O}$ be a smooth orbifold. Let
  $x\in\orbify{O}$. Denote by $A_{x}\orbify{O}$ the set of {\em
    admissible tangent vectors at $x$}
$$A_{x}\orbify{O}=\left\{v\in T_x\orbify{O}\mid (x,v)=%
  \sigma(x)\textup{ for some
  }\sigma\in\mathscr{D}^0_{\Orb}(\orbify{O})\right\}\subset
T_x\orbify{O}$$ By proposition~\ref{OrbiSecsAreVecSpace},
$A_{x}\orbify{O}$ is a vector space for each $x$, and a suborbifold of
$T_x\orbify{O}$.  The {\em admissible tangent bundle of $\orbify{O}$}
is the subset
$A\orbify{O}=\bigcup_{x\in\orbify{O}}A_x\orbify{O}\subset T\orbify{O}$
with the subspace topology.  It is not hard to see that, in general,
$A\orbify{O}$ is not an orbifold. See example~\ref{RmodZ2}.

\end{definition}

\begin{figure}[h]
  \centering
  \includegraphics[scale=.45]{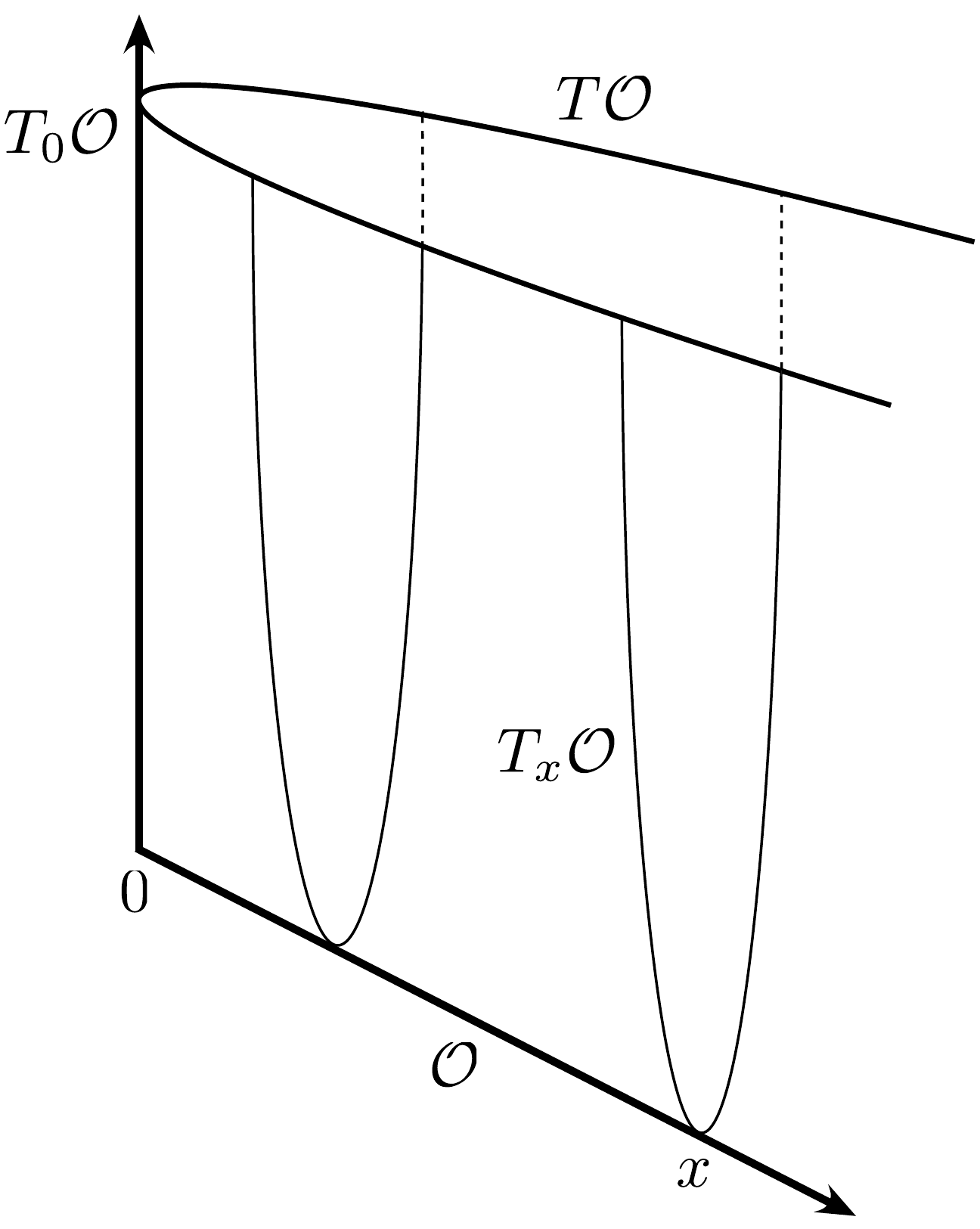}
  \hspace{.2in}
  \includegraphics[scale=.5]{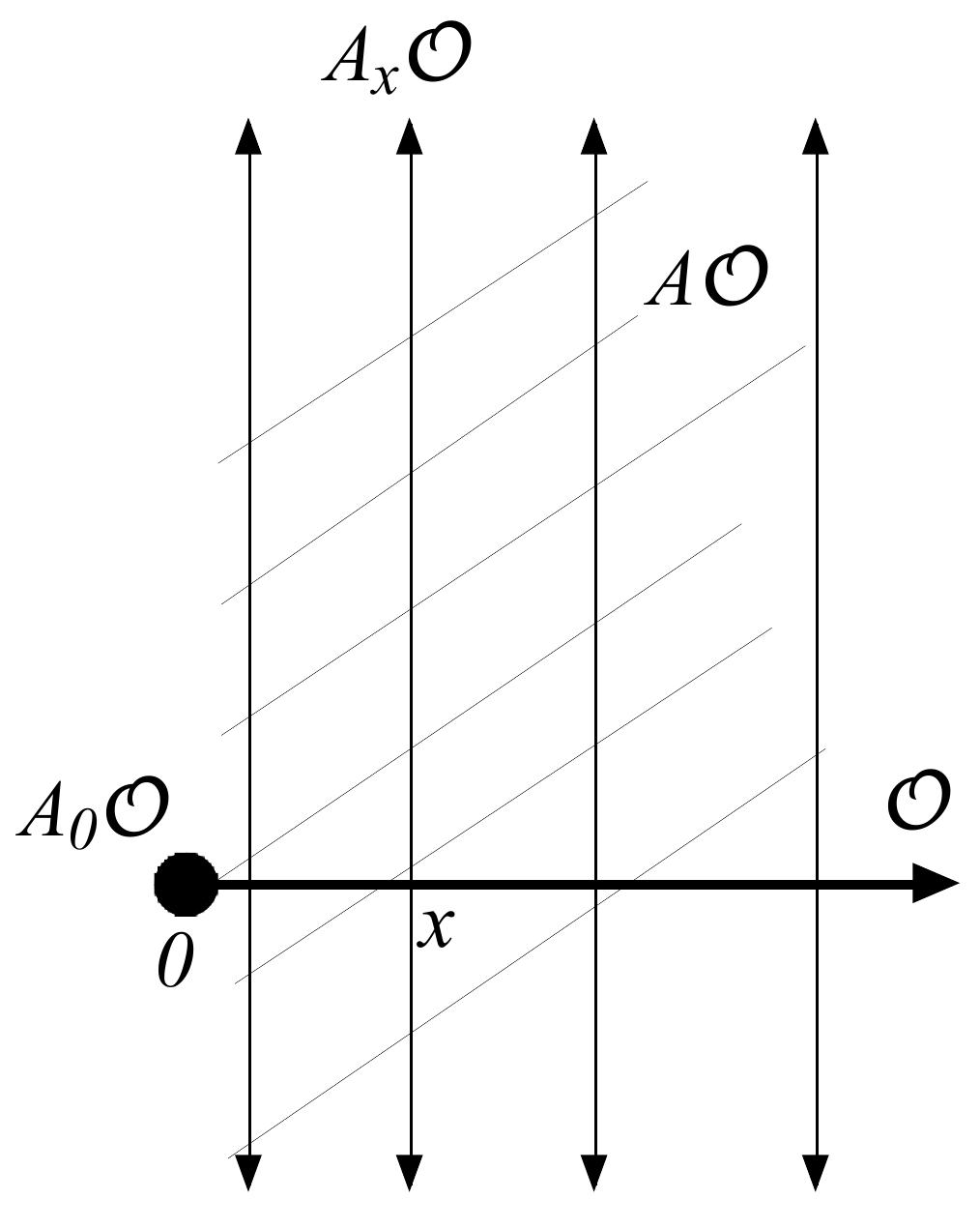}
  \caption{The tangent and admissible tangent bundles of
    example~\ref{RmodZ2}}
  \label{TangentBundleExample}
\end{figure}

\begin{example}\label{RmodZ2} Let $\orbify{O}$ be the orbifold
  $\R/\Z_2$ where $\Z_2$ acts on $\R$ via $x \to -x$. The underlying
  topological space $X_{\orbify{O}}$ of $\orbify{O}$ is $[0, \infty)$
  and the isotropy subgoups are trivial for $x\in (0, \infty)$ and
  $\Z_2$ for $x = 0$.  The tangent orbibundle $T\orbify{O}$ is given
  by $(\R\times\R)/\Z_2$ with the $\Z_2$ action being given by $(x,
  y)\to(-x, -y)$, with underlying topological space the quotient of
  $[0, \infty)\times\R$ by the equivalence relation $(0, y) \sim (0,
  -y)$. Note that $T_x\orbify{O} = \R$ if $x\ne 0$ but that
  $T_0\orbify{O} = [0, \infty)$. It also follows from
  proposition~\ref{OrbiSecsAreVecSpace} that the set of admissible
  vectors at $x=0$ consists only of the zero vector.  Thus, all
  orbisections $\sigma\in\mathscr{D}^r_{\Orb}(\orbify{O})$ must vanish
  at $0$. In particular, $A\orbify{O}\cong \{(0,0)\}\cup
  \{(0,\infty)\times\R\}$ and a neighborhood of $(0,0)$ is not covered
  by an orbifold chart, and thus $A\orbify{O}$ is not an orbifold.
  See Figure~\ref{TangentBundleExample}.
\end{example}

\begin{proposition} For a compact orbifold $\orbify{O}$, the inclusion
  $\mathscr{D}^r_{\Orb}(\orbify{O})\hookrightarrow
  C^r_{\Orb}(\orbify{O}, T\orbify{O})$ induces a separable Banach
  space structure on $\mathscr{D}^r_{\Orb}(\orbify{O})$ for $1\le
  r<\infty$ and a separable \Frechet space structure if $r=\infty$.
\end{proposition}
\begin{proof}

  Let $\mathcal{C} = \{C_i\}_{i = 1}^N$ be a cover of $\orbify{O}$ by
  a finite number of compact charts (obtained by passing to a finite
  subcover of a covering by orbifold charts and then shrinking if
  necessary), equipped with trivializations $\Psi_i:T_{
    C_i}\orbify{O}\to (\tilde C_i\times\R^n)/\Gamma_i$ of the tangent
  orbibundle over $ C_i$ where the lifts $\tilde\Psi_i$ are linear in
  the fiber.  Let $V_{i,r} = C^r(\tilde C_i,\R^n)$ for $i = 1,\ldots,
  N$ and $0\le r\le\infty$ with topology of uniform convergence of
  derivatives of order $\le r$.  This is a Banach space for finite $r$
  and a \Frechet space for $r=\infty$. For finite $r$, let $\|\
  \|_{i,r}$ be a $C^r$ norm on $V_{i,r}$.  Define a linear map
  $T:\mathscr{D}^r_{\Orb}(\orbify{O})\to\bigoplus_{i = 1}^N V_{i,r}$
  by
  \begin{equation*}
    T(\sigma) = \left(\text{pr}_2(\tilde\Psi_1(\tilde\chi_1\tilde\sigma)), \ldots, \text{pr}_2(\tilde\Psi_N(\tilde\chi_N\tilde\sigma))\right)
  \end{equation*}
  where $\chi_i\in C_{\Orb}^r(\orbify{O},[0,1])$, $i = 1,\ldots,N$, is
  a partition of unity subordinate to the cover $\mathcal{C}$ (see
  proposition~\ref{PartitionsOfUnity} for a proof of the existence of
  such partitions of unity) and $\text{pr}_2:\tilde
  C_i\times\R^n\to\R^n$ is bundle projection onto the second factor.
  Continuity of $T$ is immediate from the definitions of the $C^r$
  topology on $\mathscr{D}^r_{\Orb}(\orbify{O})$ and the topology on
  $\bigoplus_{i =1}^N V_{i,r}$.  Moreover, given a neighborhood of the
  zero section $\mathbf{0}\in\mathscr{D}^r_{\Orb}(\orbify{O})$ of the
  form $\mathscr{N}^r(\mathbf{0},\varepsilon_i;\mathcal{C})$, it is
  apparent that there is a neighborhood of the zero section
  $\mathbf{0}$ in $\bigoplus_{i=1}^N V_{i,r}$ of the form
  $\max\{\|s_1\|_{1,r},\ldots,\|s_N\|_{N,r}\} < \delta$ where $\delta
  \le \min\{ \varepsilon_1,\ldots,\varepsilon_N\}$ contained in
  $T\left(\mathscr{N}^r(\mathbf{0},\varepsilon_i;\mathcal{C})\right)$.
  Thus, with the subspace topology on
  $T(\mathscr{D}^r_{\Orb}(\orbify{O}))$,
  $T:\mathscr{D}^r_{\Orb}(\orbify{O})\to
  T(\mathscr{D}^r_{\Orb}(\orbify{O}))$ is a linear homeomorphism.
  Since $\mathscr{D}^r_{\Orb}(\orbify{O})\subset
  C^r_{\Orb}(\orbify{O},T\orbify{O})$ is a closed subset, we see that
  $T(\mathscr{D}^r_{\Orb}(\orbify{O}))$ is a closed subspace of the
  direct sum and thus $\mathscr{D}^r_{\Orb}(\orbify{O})$ inherits a
  Banach space structure if $r <\infty$ and a \Frechet space structure
  if $r=\infty$.
\end{proof}

\subsection*{Curves in Orbifolds}\label{CurvesInOrbifoldsSection}

In this paragraph we study the notion of curves in orbifolds. As a
special case of example~\ref{MapsFromManifoldsToOrbifolds} we make the
following

\begin{definition}\label{OrbifoldCurves} Let $I$ be an interval
  (finite or infinite, closed, open or half-open) with trivial
  orbifold structure and $\orbify{O}$ a smooth orbifold.  Then
  elements of $C^r_{\Orb}(I,\orbify{O})$ are the {\em $C^r$ orbifold
    curves in $\orbify{O}$}.
\end{definition}

\begin{definition}\label{TangentVectorToCurve} Let $\orbify{O}$ be a
  smooth $C^{r+1}$ orbifold, and let $c\in C^r_{\Orb}(I,\orbify{O})$
  be an orbifold curve. Suppose $\tilde c_{\tilde x}$ is a $C^r$ lift
  of $c$ to a chart $\tilde U_x$. Let $\tilde c'_{\tilde x}(t)$ be the
  tangent vector at $t$.  If $\Pi_x\left(\tilde c_{\tilde x}(t),\tilde
    c'_{\tilde x}(t)\right)=\left(c(t),v\right)\in TU_x$, then $v\in
  T_{c(t)}U_x$ is called the {\em tangent vector to $c$ at $t$} and we
  denote it by $c'(t)$.
\end{definition}

\begin{proposition}\label{TangentVectorWellDefined} If $c\in
  C^r_{\Orb}(I,\orbify{O})$, then the tangent vector $c'(t)$ is
  well-defined.
\end{proposition}
\begin{proof} Let $x_0=c(t_0)$ and consider an orbifold chart $(\tilde
  U_{x_0},\Gamma_{x_0})$ at $x_0$. Let $t_0\in J\subset I$ be an
  interval such that $c(t)\in U_{x_0}$ for all $t\in J$. Let
  $\hat{c}(t)$ be a $C^r$ lift of $c(t)$ to $\tilde U_{x_0}$. If $x_0$
  is non-singular, then $\Gamma_{x_0}$ is trivial and $\hat{c}(t)$ is
  unique. Thus, $c'(t_0)$ is well defined when $x_0$ is non-singular.

  Now suppose that $x_0$ is singular. If $t_0\in\partial I$, it is not
  hard to see (since $\Gamma_{x_0}$ is finite, acts discretely, and
  lifts are continuous) that there is a subinterval $t_0\in J'\subset
  J$ such that any other lift of $c(t)$ is of the form $\tilde
  c(t)=\gamma\cdot\hat{c}(t)$. This is a $C^r$ lift of $c$ for any
  $\gamma\in\Gamma_{x_0}$. The tangent vector $\tilde
  c'(t_0)=d\gamma_{\hat{c}(t_0)}\hat{c}'(t_0)$. Thus, $\tilde c'(t_0)$
  is in the same orbit as $\hat{c}'(t_0)$ of the $\Gamma_{x_0}$ action
  on $T_{\tilde x_0}\tilde U_{x_0}$ and so their projections to
  $T_{x_0}U_{x_0}$ are equal and thus $c'(t_0)$ is well-defined.  If
  $t_0$ is an interior point of $I$, then it is possible to build a
  $C^0$ lift of $c$ by concatenation:

$$\tilde c(t)=
\begin{cases} 
  \hat{c}(t) & \text{for $t\le t_0$}\\
  \gamma\cdot\hat{c}(t) & \text{for $t\ge t_0$}
\end{cases}
$$
Note that by our previous observations this is the only way to produce
another lift around $t_0$. The condition that $\tilde c$ be at least
$C^1$ implies that
$\hat{c}'(t_0)=d\gamma_{\hat{c}(t_0)}\hat{c}'(t_0)$. Thus, like above,
we see that $c'(t_0)$ is well defined and furthermore that
$\hat{c}'(t_0)$ is fixed by the action of $\gamma$ on $T_{\tilde
  x_0}\tilde U_{x_0}$. Note that $c'(t_0)$ is not necessarily an
admissible tangent vector, as $\hat{c}'(t_0)$ is not necessarily fixed
by all elements of $\Gamma_{x_0}$.
\end{proof}

\begin{example}\label{CurveLiftings} Let $\orbify{O}$ be the orbifold
  $\R^2/\Z_2$ where $\Z_2$ acts on $\R^2$ via $(x,y) \to (x,-y)$. The
  underlying topological space $X_{\orbify{O}}$ of $\orbify{O}$ is the
  closed upper half-plane and the isotropy subgoups at $(x,y)$ are
  $\Z_2$ if $y=0$ and trivial otherwise. Let $I=[-1,1]$ and consider
  the curves $b(t)=(t,|t|)$ and $c(t)=(t,t^2)$. It's easy to see that
  $b$ and $c$ have four $C^0$ lifts. They are of the form:

$$\tilde b^\pm_\pm(t)=
\begin{cases} (t,\pm t) & \text{for $t\le 0$}\\
  (t,\pm t) & \text{for $t\ge 0$}
\end{cases}
\hspace{.5in} \tilde c^\pm_\pm(t)=
\begin{cases} (t,\pm t^2) & \text{for $t\le 0$}\\
  (t,\pm t^2) & \text{for $t\ge 0$}
\end{cases}
$$
$b$ has two $C^r$ lifts, $b^+_+$ and $b^-_-$ for $r\ge 1$. However,
all four lifts of $c$ are $C^1$ while only two, $c^+_+$ and $c^-_-$,
are $C^r$ for $r\ge 2$. One sees that in the case of $b$ the $C^1$
lifts do not arise from a non-trivial concatenation. Note that the
tangent vectors of these lifts at $t=0$ are not fixed by the action of
$\Z_2$. On the other hand, two of the four $C^1$ lifts of $c$ do arise
as non-trivial concatenations. Their tangent vectors at $t=0$ are
fixed by the $\Z_2$ action.
\end{example}

\section{Smooth Riemannian Orbifold
  Structures}\label{SmoothRiemannianStructureSection}

In this section we show that any smooth orbifold admits a smooth
Riemannian orbifold structure.  Although orbifolds are metrizable,
this is not sufficient for our needs as we will need to make use of a
smooth orbifold Riemannian exponential map:
$\exp:T\orbify{O}\to\orbify{O}$. In order to do this, we proceed as in
the classical situation of Riemannian manifolds.

\begin{proposition}\label{PartitionsOfUnity} Let $\orbify{O}$ be a
  smooth orbifold and let $\mathscr{U}=\{U_\alpha\}_{\alpha\in I}$ be
  a locally finite open covering of $\orbify{O}$ by orbifold charts.
  Then there exists a $C^\infty$ partition of unity subordinate to
  $\mathscr{U}$.
\end{proposition}

\begin{proof}
  Paracompactness of $\orbify{O}$ implies the existence of the
  covering $\mathscr{U}$.  Without loss of generality, by
  proposition~\ref{OrbifoldsAreSmoothable}, we may assume that
  $\orbify{O}$ is a $C^\infty$ orbifold.  Let
  $\tilde{\mathscr{U}}=\{(\tilde U_\alpha,\Gamma_\alpha)\}_{\alpha\in
    I}$ be the corresponding covering charts and let
  $\pi_\alpha:\tilde U_\alpha\to U_\alpha$ be the quotient map.  Since
  $\orbify{O}$ is paracompact and Hausdorff, we let
  $\{\chi'_\alpha\}:\orbify{O}\to [0,1]$ be a $C^0$ partition of unity
  subordinate to the cover $\{U_\alpha\}$. If we give $[0,1]$ the
  trivial orbifold structure, we may regard each $\chi'_\alpha$ as an
  element of $C^0_{\Orb}(\orbify{O},[0,1])$ (See
  example~\ref{TopMapsFromOrbifoldsToManifolds}). That is, each local
  lift of $\chi_\alpha$, $\tilde\chi'_{\alpha,\beta}:\tilde U_\beta\to
  [0,1]$, is $C^0$ equivariant and $\tilde\chi'_{\alpha,\beta}(\tilde
  x)=\chi'_\alpha\circ\pi_\beta(\tilde x)$ for all $\tilde x\in\tilde
  U_\beta$. Note that for fixed $x\in\orbify{O}$,
  $\pi^{-1}_\beta(x)\ne\emptyset$ for only finitely many $\beta$ and
  furthermore,
  $\tilde\chi'_{\alpha,\beta}\left(\pi^{-1}_\beta(x)\right)\ne 0$ for
  all but a finite number of $\alpha$.  In order to produce a
  $C^\infty$ partition of unity we choose, for each pair
  $(\alpha,\beta)$, a nonnegative $C^\infty$ map
  $\tilde\chi''_{\alpha,\beta}:\tilde U_\beta\to [0,1]$ which is
  sufficiently $C^0$ close to $\tilde\chi'_{\alpha,\beta}$. For
  $\tilde x\in\tilde U_\beta$ define
  $$\hat\chi_{\alpha,\beta}(\tilde x)=%
  \frac{1}{|\Gamma_\beta|}\sum_{\gamma\in\Gamma_\beta}\tilde\chi''_{\alpha,\beta}(\gamma\cdot\tilde
  x)$$ By defining
  $\tilde\chi_{\alpha,\beta}=\displaystyle\frac{\hat\chi_{\alpha,\beta}}{\sum_{\mu,\nu\in
      I}\hat\chi_{\mu,\nu}}$ we get a $C^\infty$
  $\Gamma_\beta$-equivariant map on $\tilde U_\beta$ that is $C^0$
  close to $\tilde\chi'_{\alpha,\beta}$ for each pair
  $(\alpha,\beta)$.  Thus the map
  $$\chi_\alpha(x)=
    \begin{cases}
      \sum_\beta\tilde\chi_{\alpha,\beta}\left(\pi_\beta^{-1}(x)\right) & \text{for $x\in U_\alpha$}\\
      0 & \text{for $x\in\orbify{O}-U_\alpha$}
    \end{cases}$$ is well-defined, each $\chi_\alpha\in
    C^\infty_{\Orb}(\orbify{O},[0,1])$ and the collection
    $\{\chi_\alpha\}$ is a smooth partition of unity subordinate to
    the cover $\{U_\alpha\}$.
  \end{proof}

  We now prove the existence of a smooth orbifold Riemannian metric.
  We could, of course, do this by defining appropriate notions of
  tensor bundles over orbifolds and their sections. However, since our
  needs are limited, we choose to do this in an elementary way
  following the classical development. Since the tangent space
  $T_x\orbify{O}\cong\R^n/\Gamma_x$ is, in general, a convex cone
  rather than a vector space, we make the following

\begin{definition}\label{OrbifoldInnerProduct}
  A function $g_x:T_x\orbify{O}\times T_x\orbify{O}\to\R$ is a {\em
    positive definite, real, orbifold inner product} if it has a
  $\Gamma_x\times\Gamma_x$ equivariant lift $\tilde
  g_x:\R^n\times\R^n\to\R$ which is a positive definite real inner
  product on $\R^n$.  Note that we gave the natural product orbifold
  structure to $T_x\orbify{O}\times T_x\orbify{O}$.
\end{definition}

\begin{definition}\label{OrbifoldRiemannianMetric} Let $\orbify{O}$ be
  a smooth $C^{r+1}$ orbifold.  A {\em smooth $C^r$ orbifold
    Riemannian metric on} $\orbify{O}$ is a collection
  $g=\{g_x\}_{x\in\orbify{O}}$ of positive definite real orbifold
  inner products so that the functions $g(\sigma,\tau):x\mapsto
  g_x(\sigma(x),\tau(x))$ are elements of $C^r_{\Orb}(\orbify{O},\R)$
  for all orbisections
  $\sigma,\tau\in\mathscr{D}^r_{\Orb}(\orbify{O})$. An orbifold
  equipped with a $C^r$ Riemannian metric will be called a {\em $C^r$
    Riemannian orbifold}.
\end{definition}

\begin{proposition}\label{OrbifoldRiemannianMetricExists}
  Let $\orbify{O}$ be a smooth orbifold. Then there exists on
  $\orbify{O}$ a smooth $C^\infty$ orbifold Riemannian metric.
\end{proposition}

\begin{proof} Without loss of generality, by
  proposition~\ref{OrbifoldsAreSmoothable}, we may assume that
  $\orbify{O}$ is a $C^\infty$ orbifold.  Using the notation from
  proposition~\ref{PartitionsOfUnity}, let $\{\chi_\alpha\}$ be a
  $C^\infty$ partition of unity and let $\tilde g'_\alpha$ be a
  $C^\infty$ Riemannian metric on $\tilde U_\alpha$. Define
$$\tilde g_\alpha(\tilde v,\tilde w)=\frac{1}{|\Gamma_\alpha|^2}\sum_{(\gamma,\mu)\in\Gamma_\alpha\times\Gamma_\alpha}%
\tilde g'_\alpha(d\gamma_{\tilde x}(\tilde v),d\mu_{\tilde x}(\tilde
w))$$ for $\tilde v,\tilde w\in T_{\tilde {x}}\tilde U_\alpha$.  Then
$\tilde g_\alpha$ is a $C^\infty$, $\Gamma_\alpha\times\Gamma_\alpha$
equivariant positive definite, real inner product on each $T_{\tilde
  x}\tilde U_\alpha$ which descends to a smooth orbifold Riemannian
metric $g_\alpha$ on $U_\alpha$. Thus, $g=\sum_\alpha\chi_\alpha
g_\alpha$ is the required $C^\infty$ orbifold Riemannian metric on
$\orbify{O}$.
\end{proof}

\begin{remark}\label{RiemanianOrbifoldDefisConsistent} Note that the
  proof of proposition~\ref{OrbifoldRiemannianMetricExists} shows that
  the action of $\Gamma_\alpha$ on $\tilde U_\alpha$ is by isometries
  relative to $\tilde g_\alpha$, and that the equivariant transition
  maps $\tilde\psi$ that appear in definition~\ref{orbifold} are
  isometric embeddings.  By shrinking the cover $\{U_\alpha\}$ if
  necessary, we may assume that each orbifold covering chart $\tilde
  U_\alpha$ is convex making $\orbify{O}$ a Riemannian orbifold as
  defined in \cite{MR1218706} and \cite{BorPhD}. Recall that for a
  Riemannian manifold to be convex means there exists a unique minimal
  geodesic joining any two points.
\end{remark}

If $\orbify{O}$ is a smooth $C^r$ Riemannian orbifold, then we may
give $\orbify{O}$ the structure of a length space.  A general
reference is \cite{MR1699320}.  In particular, given two points
$x,y\in\orbify{O}$ we may define the distance between $x$ and $y$ to
be
$$d(x,y)=\inf\{\textup{Length}(c)\mid c\in C^0_{\Orb}(I,\orbify{O})\textup{ and } c \textup{ joins } x\textup{ to } y\}$$
The length of a curve $c$ is defined by adding up the lengths of local
lifts in each orbifold chart $\tilde U_\alpha$.  This can be shown to
be well-defined and independent of the choice of lift \cite{BorPhD}.
This length metric structure generates a topology that is the same as
the as the topology of the underlying space of $\orbify{O}$. If
$(\orbify{O},d)$ is complete any two points can be joined by a {\em
  minimal geodesic} realizing the distance $d(x,y)$ \cite{MR1699320},
since $\orbify{O}$ is locally compact.  Moreover, the local lifts of
any such minimal geodesic must be a smooth $C^r$ minimal geodesic in
each $\tilde U_\alpha$, justifying the use of the terminology.
Additionally, if $c\in C^r_{\Orb}(I,\orbify{O})$ is a minimal geodesic it
can be shown that $\Gamma_{c(t)}=\Gamma_{c(t')}$ for all $t,t'\in
I-\partial I$ \cite{MR1218706}.

We now proceed to define the exponential map for a Riemannian
orbifold.  For a general reference for standard results of Riemannian
geometry that we need see \cite{MR1480173}.  As in the proof of
proposition~\ref{OrbifoldRiemannianMetricExists}, assume the
collection $\{U_\alpha\}$ is a locally finite open covering of
$\orbify{O}$ by orbifold charts that are relatively compact.  Let
$TU_\alpha \cong (\tilde U_\alpha\times\R^n)/\Gamma_\alpha$ be a local
trivialization of the tangent bundle over $U_\alpha$.  Denote the
Riemannian exponential map on $T\tilde U_\alpha$ by
$\widetilde{\exp}_{\tilde U_\alpha}:T\tilde U_\alpha\to \tilde
U_\alpha$.  Thus, for $\tilde x\in\tilde U_\alpha$ and $\tilde v\in
T_{\tilde x}\tilde U_\alpha$ we have $\widetilde{\exp}_{\tilde
  U_\alpha}(\tilde x,t\tilde v)=\tilde c_{\tilde x,\tilde v}(t)$ where
$\tilde c_{\tilde x,\tilde v}$ is the unit speed geodesic in $\tilde
U_\alpha$ which starts at $\tilde x$ and has initial velocity $\tilde
v$.  Recall that there is an open neighborhood $\tilde\Omega_{\tilde
  U_\alpha}\subset T\tilde U_\alpha$ of the $0$-section of $T\tilde
U_\alpha$ such that $\tilde c_{\tilde x,\tilde v}(1)$ is defined for
$\tilde v\in T_{\tilde x}\tilde U_\alpha\cap\tilde\Omega_{\tilde
  U_\alpha}$. Furthermore, by shrinking $\tilde\Omega_{\tilde
  U_\alpha}$ if necessary, we may assume that on $T_{\tilde x}\tilde
U_\alpha\cap\tilde\Omega_{\tilde U_\alpha}$, $\widetilde{\exp}_{\tilde
  U_\alpha}(\tilde x,\cdot)$ is a local diffeomorphism onto a
neighborhood of $\tilde x\in\tilde U_\alpha$ for each $\tilde
x\in\tilde U_\alpha$.  Let
$\Omega_\alpha=\Pi_\alpha(\tilde\Omega_{\tilde U_\alpha})$, an open
subset of $T\orbify{O}$, and define
$\Omega=\bigcup_{\alpha}\Omega_\alpha$. $\Omega$ is an open
neighborhood of the $0$-orbisection of $T\orbify{O}$.

\begin{definition}\label{ExponentialMap} Let $x\in U_\alpha$, and
  $(x,v)\in\Omega_\alpha$.  Choose $(\tilde x,\tilde
  v)\in\Pi_\alpha^{-1}(x,v)$. Then the {\em Riemannian exponential
    map} $\exp:\Omega\subset T\orbify{O}\to\orbify{O}$ is defined by
  $\exp(x,v)=\pi_\alpha\circ\widetilde{\exp}_{\tilde U_\alpha}(\tilde
  x,\tilde v)$.
\end{definition}

\begin{proposition}\label{ExpMapWellDefined} Let $\orbify{O}$ be a
  $C^{r+1}$ Riemannian orbifold. Then the exponential map
  $\exp(x,v)=\pi_\alpha\circ\widetilde{\exp}_{\tilde
    U_\alpha}\circ\Pi_\alpha^{-1}(x,v)$ is well--defined.
\end{proposition}

\begin{proof} Since the metric $\tilde g_\alpha$ is equivariant
  relative to the action of $\Gamma_\alpha$ by isometries on $\tilde
  U_\alpha$ we see that (since isometries map geodesics to geodesics)
  $\widetilde{\exp}_{\tilde U_\alpha}[\gamma\cdot(\tilde x,\tilde
  v)]=\gamma\cdot\widetilde{\exp}_{\tilde U_\alpha}(\tilde x,\tilde
  v)$.  Thus, $\widetilde{\exp}_{\tilde
    U_\alpha}:\tilde\Omega_\alpha\subset T\tilde U_\alpha\to\tilde
  U_\alpha$ is equivariant and hence $\exp$ is well-defined for each
  $\tilde U_\alpha$.  If $x\in U_\alpha\cap U_\beta$, then there is an
  orbifold chart $U_{\alpha\beta}\subset U_\alpha\cap U_\beta$ of $x$,
  and equivariant isometric embeddings $\tilde\psi_{\alpha}:\tilde
  U_{\alpha\beta}\to\tilde U_\alpha$ and $\tilde\psi_{\beta}:\tilde
  U_{\alpha\beta}\to\tilde U_\beta$. This observation is enough to
  show that $\exp$ is independent of local chart.
\end{proof}

As usual we denote by $\exp_x$ the restriction of $\exp$ to a single
tangent cone $T_x\orbify{O}$. We let $B(x,r)$ denote the metric
$r$-ball centered at $x$ and use tildes to denote corresponding points
in local coverings.

\begin{proposition}\label{ExpMapHomeo} Let $\orbify{O}$ be a $C^{r+1}$
  Riemannian orbifold. Then $\exp_x$ is a local (topological)
  homeomorphism.  That is, there exists $\varepsilon>0$ such that
  $\exp_x:B(0,\varepsilon)\subset T_x\orbify{O}\to
  B(x,\varepsilon)\subset\orbify{O}$ is a (topological) homeomorphism
  with $C^r$ local lifts for each $x\in\orbify{O}$.
\end{proposition}

\begin{proof} First note that a lift of $\exp_x$ to $\tilde U_x$ is of
  the form $\widetilde{\exp}_{\tilde U_x}(\tilde x,\cdot)$.  Since the
  classical Riemannian exponential map is as smooth as its tangent
  bundle, we see that $\exp_x$ has local $C^r$ lifts.

  Choose $\varepsilon>0$ so that $\overline{B(\tilde
    0,\varepsilon)}\subset\tilde\Omega_{\tilde U_x}\cap T_{\tilde
    x}\tilde U_x$. Then $\widetilde{\exp}_{\tilde U_\alpha}(\tilde
  x,\cdot)$ is a local $C^r$ diffeomorphism from $B(\tilde
  0,\varepsilon)\subset T_{\tilde x}\tilde U_x$ onto $B(\tilde
  x,\varepsilon)\subset\tilde U_x$. By construction of the length
  metric on $\orbify{O}$, it is easy to see that
  $\pi_x\left(\overline{B(\tilde
      x,\varepsilon)}\right)=\overline{B(x,\varepsilon)}$, thus
  $\exp_x$ maps $\overline{B(0,\varepsilon)}\subset T_x\orbify{O}$
  onto $\overline{B(x,\varepsilon)}\subset\orbify{O}$.

  To see that $\exp_x$ is injective, suppose that
  $\exp_x(v)=\exp_x(w)$ for $v,w\in \overline{B(0,\varepsilon)}$. Then
  there is $\gamma\in\Gamma_x$ such that $\widetilde{\exp}_{\tilde
    U_x}(\tilde x,\tilde v)=\gamma\cdot\widetilde{\exp}_{\tilde
    U_x}(\tilde x,\tilde w)=%
  \widetilde{\exp}_{\tilde U_x}(\gamma\cdot\tilde x,d\gamma_{\tilde
    x}\tilde w)=\widetilde{\exp}_{\tilde U_x}(\tilde x,d\gamma_{\tilde
    x}\tilde w)$.  Thus, $\tilde v=d\gamma_{\tilde x}\tilde w$, since
  $\widetilde{\exp}_{\tilde U_x}(\tilde x,\cdot)$ is a local
  diffeomorphism and therefore $v=w$.

  Finally since $\exp_x$ is continuous, bijective and $\overline
  {B(0,\varepsilon)}$ is compact, we see that $\exp_x$ is a local
  homeomorphism.
\end{proof}

If we restrict the exponential map $\exp_x$ to admissible vectors at
$x$, we can say a little more.

\begin{proposition}\label{ExpAdmissibleSameStrata} Let $\orbify{O}$ be
  a $C^{r+1}$ Riemannian orbifold. Let $\varepsilon >0$ be as in
  proposition~\ref{ExpMapHomeo}.  Then the restriction of $\exp_x$ to
  $B(0,\varepsilon)\cap A_x\orbify{O}$ is a $C^r$ local diffeomorphism
  of $A_x\orbify{O}$ (with trivial suborbifold structure) onto a
  neighborhood of $x$ in the stratum $\orbify{S}_x$ (with trivial
  suborbifold structure).
\end{proposition}

\begin{proof} Let $v\in B(0,\varepsilon)\cap A_x\orbify{O}$, and
  choose $(\tilde x,\tilde v)\in\Pi_x^{-1}(x,v)\cap B(\tilde
  0,\varepsilon)$. Then, by the proof of
  proposition~\ref{OrbiSecsAreVecSpace}, $d\gamma_{\tilde x}\tilde
  v=\tilde v$ for all $\gamma\in\Gamma_x$. Thus, by equivariance of
  $\widetilde{\exp}_{\tilde U_x}$, we have for $t\in[0,1]$,
$$\widetilde{\exp}_{\tilde U_x}(\tilde x,t\tilde v)=\widetilde{\exp}_{\tilde U_x}[\gamma\cdot(\tilde x,t\tilde v)]=%
\gamma\cdot\widetilde{\exp}_{\tilde U_x}(\tilde x,t\tilde v)$$ Hence,
$\widetilde{\exp}_{\tilde U_x}(\tilde x,t\tilde v)$ is fixed by the
action of $\Gamma_x$ for all $t\in [0,1]$. This implies that for all
$t\in[0,1]$ we have,
$\exp_x(tv)=\exp(x,tv)=\pi_x\circ\widetilde{\exp}_{\tilde U_x}(\tilde
x,t\tilde v)\in B(x,\varepsilon)\cap\orbify{S}_x$.  Thus, $\exp_x$
maps onto $B(x,\varepsilon)\cap\orbify{S}_x$. In fact, since the
restriction of the action of $\Gamma_x$ to $\tilde{\orbify{S}}_x$ is
trivial ($\Gamma_x\cdot\tilde s=\tilde s$ for all $\tilde
s\in\tilde{\orbify{S}}_x$), we may identify
$\orbify{S}_x\subset\orbify{O}$ with
$\tilde{\orbify{S}}_x\subset\tilde U_x$ and under this identification
our restriction of $\exp_x$ to $A_x\orbify{O}$ is nothing more than
the map $\widetilde{\exp}_{\tilde U_x}(\tilde x,\cdot)$ restricted to
$T_{\tilde x}\tilde{\orbify{S}}_x\cap T_{\tilde x}\tilde U_x$. Hence
$\exp_x$ is a local $C^r$ (manifold) diffeomorphism.
\end{proof}

The composition of the exponential map with an orbisection turns out
to be a smooth orbifold map.

\begin{proposition}\label{ExpOrbisectionIsMap}
  Let $\orbify{O}$ be a $C^{r+1}$ Riemannian orbifold. Let $\sigma$ be
  a $C^r$ orbisection of $T\orbify{O}$. Then the map
  $E^\sigma(x)=(\exp\circ\,\sigma)(x):\orbify{O}\to\orbify{O}$ is a
  smooth $C^r$ orbifold map, provided $\sigma(x)\in\Omega$.  That is,
  $E^\sigma\in C^r_{\Orb}(\orbify{O})$.
\end{proposition}

\begin{proof} Let $(\tilde U_x,\Gamma_x)$ be an orbifold chart at
  $x\in\orbify{O}$.  For $y\in U_x$, $\sigma(y)=(y,s(y))$ where
  $s(y)\in A_y\orbify{O}$. Then as in the proof of
  proposition~\ref{OrbiSecsAreVecSpace}, if $\tilde\sigma_x$ is a lift
  of $\sigma$, then
  $\Theta_{\sigma,x}(\delta)=\delta$ for all
  $\delta\in\Gamma_x$ and $\tilde\sigma_x(\tilde y)=(\tilde
  y,\tilde s(\tilde y))$, where $\tilde s:\tilde U_x\to\R^n$ satisfies
  $\tilde s(\delta\cdot\tilde
  y)=(d\delta)_{\tilde y}\tilde s(\tilde y)$.

  The map $\tilde E^\sigma_x=\widetilde{\exp}_{\tilde
    U_x}\circ\tilde\sigma_x$ is a $C^r$ lift of $E^\sigma$ and thus we
  need to check equivariance:
  \begin{align*}
    \tilde E^\sigma_x(\delta\cdot\tilde y) & = \widetilde{\exp}_{\tilde U_x}\big(\delta\cdot\tilde y,\tilde s(\delta\cdot\tilde y)\big)\\
    & = \widetilde{\exp}_{\tilde U_x}\big(\delta\cdot\tilde y,(d\delta)_{\tilde y}{\tilde s}(\tilde y)\big)\\
    & = \widetilde{\exp}_{\tilde U_x}\big[\delta\cdot\big(\tilde y,\tilde s(\tilde y)\big)\big]\\
    & = \delta\cdot\widetilde{\exp}_{\tilde U_x}\big(\tilde y,\tilde s(\tilde y)\big)\\
    & = \delta\cdot \tilde E^\sigma_x(\tilde y)
  \end{align*}

  Thus, $\tilde E^\sigma_x$ is $\Theta_{E^\sigma,x}$ equivariant if we
  define $\Theta_{E^\sigma,x}(\delta)=\delta$. Hence $E^\sigma\in C^r_{\Orb}(\orbify{O})$.
\end{proof}

Denote by $\mathbf{0}:\orbify{O}\to T\orbify{O}$,
$\mathbf{0}(x)=0_x\in T_x\orbify{O}$, the $0$-orbisection of
$T\orbify{O}$. The next proposition shows that if $\sigma$ is
sufficiently $C^1$ close to the $0$-orbisection $\mathbf{0}$, then
$E^\sigma$ is a local orbifold diffeomorphism.

\begin{proposition}\label{ExpOrbisectionIsLocalDiffeo}
  Let $\orbify{O}$ be a $C^{r+1}$ Riemannian orbifold and
  $U_\alpha\subset\orbify{O}$, where $U_\alpha$ is a relatively
  compact orbifold chart.  Then there is a open neighboorhood
  $\Lambda_\alpha\subset\Omega_\alpha\subset TU_\alpha$ of
  $U_\alpha\times \{0\}\subset TU_\alpha$, such that if $\sigma$ is a
  $C^r$ orbisection with $\sigma(x)\in\Lambda_\alpha$ and $\sigma$ is
  sufficiently $C^1$ close to $\mathbf{0}$ on $U_\alpha$, then
  $E^\sigma|_{U_\alpha}$ is a $C^r$ orbifold diffeomorphism onto its
  image. That is, $E^\sigma|_{U_\alpha}$ is an orbifold embedding.
\end{proposition}

\begin{proof}
  Without loss of generality, we may assume by shrinking $U_\alpha$
  and $\Lambda_\alpha$ if necessary that $U_\alpha$ and
  $E^\sigma(U_\alpha)$ are contained in a single relatively compact
  orbifold chart $(\tilde U,\Gamma)$.  Let
  $\Lambda_\alpha=\Pi_\alpha\left(\tilde\Omega_{\tilde
      U}\cap\tilde\Omega_{\tilde U_\alpha}\right)$.  By
  proposition~\ref{ExpOrbisectionIsMap}, we know already that
  $E^\sigma(x)$ is a $C^r$ orbifold map. We need to show that
  $E^\sigma$ has an inverse that is also a $C^r$ orbifold map. We
  first show that $E^\sigma(x)$ is injective.

  There exists $\gamma\in\Gamma$ such that
  $\left(\widetilde{\exp}_{\tilde
      U_u}\circ\tilde{\mathbf{0}}_u\right)(\tilde x)=\gamma\cdot\tilde
  x$ since this map is a lift of the identity map. If $\sigma$ is
  $C^1$ close enough to $\mathbf{0}$ with lift
  $\tilde{\mathbf{0}}_u=(\tilde x,0)$, then $\tilde\sigma_u=(\tilde
  x,\tilde s(\tilde x))$ for $u\in U$.  Suppose that
  $E^\sigma(x)=E^\sigma(y)=u$ for $x,y,u\in U$.  (This implies that
  the isotropy groups of $x,y,u$ are equal, by
  proposition~\ref{ExpAdmissibleSameStrata}).  Then there exists
  $\delta\in\Gamma$ such that $\tilde E^\sigma_u(\tilde
  x)=\delta\cdot\tilde E^\sigma_u(\tilde y)$.  Thus,
  \begin{alignat*}{2}
    \left(\widetilde{\exp}_{\tilde
        U_u}\circ\tilde\sigma_u\right)(\tilde x) & =%
    \delta\cdot\left[\left(\widetilde{\exp}_{\tilde U_u}\circ\tilde\sigma_u\right)(\tilde y)\right] & \\
    & =\delta\cdot\left(\widetilde{\exp}_{\tilde U_u}\big(\tilde y,\tilde s(\tilde y)\big)\right) & \\
    & =\widetilde{\exp}_{\tilde U_u}\big(\delta\cdot\tilde y,(d\delta)_{\tilde y}\tilde s(\tilde y)\big) & \\
    & =\widetilde{\exp}_{\tilde U_u}\big(\tilde w,\tilde s(\tilde w)\big) & \text{where }\tilde w=\delta\cdot\tilde y\\
    & = \left(\widetilde{\exp}_{\tilde U_u}\circ\tilde\sigma_u\right)(\tilde w) &\\
  \end{alignat*}
  Since a sufficiently small $C^1$ neighborhood of an embedding is an
  embedding \cite{MR0198479}, by choosing $\sigma$ sufficiently $C^1$
  close to $\mathbf{0}$, we may conclude that $\tilde x=\tilde w$
  which in turn implies that $\tilde x$ and $\tilde y$ are in the same
  orbit of the $\Gamma$ action on $\tilde U$. Thus $x=y$.

  We now show that $(E^\sigma)^{-1}$ is a $C^r$ orbifold map. Denote
  by $\widetilde{\exp}^{-1}_{\tilde U_u,\tilde x}$ the $C^r$ map
  $\left[\widetilde{\exp}_{\tilde U_u}(\tilde
    x,\cdot)\right]^{-1}:\tilde U\to T_{\tilde x}\tilde U$. Also, let
  $\text{pr}_1:T\tilde U\to\tilde U$ be the bundle projection $(\tilde
  x,\tilde v)\mapsto\tilde x$. Suppose $\tilde y=\tilde
  E^\sigma_u(\tilde x)$. We claim that $\left(\tilde
    E^\sigma_u\right)^{-1}(\tilde y)=\text{
  pr}_1\left(\widetilde{\exp}^{-1}_{\tilde U_u,\gamma\tilde x}(\tilde
    y)\right)$, a composition of $C^r$ maps. To see the formula is
  correct we compute:
  \begin{align*}
    \text{pr}_1\left(\widetilde{\exp}^{-1}_{\tilde U_u,\gamma\tilde
        x}(\tilde y)\right) & =%
    \text{pr}_1\left(\widetilde{\exp}^{-1}_{\tilde U_u,\gamma\tilde x}\left(\tilde E^\sigma_u(\tilde x)\right)\right)\\
    & = \text{pr}_1\left(\widetilde{\exp}^{-1}_{\tilde
        U_u,\gamma\tilde x}\left(\widetilde{\exp}_{\tilde U_u}%
        \circ\tilde\sigma_u(\tilde x)\right)\right)\\
    & = \text{pr}_1\left(\widetilde{\exp}^{-1}_{\tilde
        U_u,\gamma\tilde x}\left(\widetilde{\exp}_{\tilde U_u}%
        (\tilde x,\tilde s(\tilde x)\right)\right)\\
    & = \text{pr}_1\left(\tilde x,\tilde s(\tilde x)\right)\\
    & =\tilde x
  \end{align*}

  Now we need to check equivariance. From the computation in
  proposition~\ref{ExpOrbisectionIsMap}, for any $\delta\in\Gamma$, we
  have $\tilde E^\sigma_u\left(\delta\cdot\tilde
    x\right)=\delta\cdot\tilde y$. Thus,
  \begin{align*}
    \left(\tilde E^\sigma_u\right)^{-1}(\delta\cdot\tilde y) & =
    \text{pr}_1\left(\widetilde{\exp}^{-1}_{\tilde
        U_u,\delta\gamma\tilde x}%
      \left(\tilde E^\sigma_u\left(\delta\cdot\tilde x\right)\right)\right)\\
    & = \text{pr}_1\left(\widetilde{\exp}^{-1}_{\tilde
        U_u,\delta\gamma\tilde x}%
      \left[\widetilde{\exp}_{\tilde U_u}\left(\delta\cdot\tilde x,\tilde s\left(\delta\cdot\tilde x\right)\right)\right]\right)\\
    & = \text{pr}_1\left(\left(\delta\cdot\tilde x,\tilde s\left(\delta\cdot\tilde x\right)\right)\right)\\
    & = \delta\cdot\tilde x\\
    & = \delta\cdot\left(\tilde E^\sigma_u\right)^{-1}(\tilde y)
  \end{align*}

  Thus, $\left(\tilde E^\sigma_u\right)^{-1}$ is
  $\Theta_{(E^\sigma)^{-1},u}$ equivariant if we define
  $\Theta_{(E^\sigma)^{-1},u}(\delta)=\delta$. Note that
  $\Theta_{(E^\sigma)^{-1},u}=\left(\Theta_{E^\sigma,u}\right)^{-1}$
  as to be expected.
\end{proof}

The next lemma is a standard result of differential topology adapted to orbifolds:

\begin{lemma}\label{C0CloseToIdIsOnto}
  Let $\textup{Id}:\orbify{O}\to\orbify{O}$ be the identity map. Then
  there is a $C^0$ neighborhood of $\textup{Id}$ such that if $f$ lies
  in this neighborhood, then $f$ is surjective.
\end{lemma}
\begin{proof} The proof is essentially a minor modification of the
  argument in \cite[lemma 3.11]{MR0198479}. For completeness, we give
  it here.  Let $\{C_i\}$ be a locally finite covering of $\orbify{O}$
  by compact sets whose interiors also cover $\orbify{O}$. Assume
  further that the corresponding orbifold charts $(\tilde
  C_i,\Gamma_i)$ have $\tilde C_i=\text{unit ball }B^n\subset\R^n$,
  and let $(\tilde V_i,\Gamma_i)$ be an orbifold chart with $\tilde
  C_i\subset\intr(\tilde V_i)$. Let $\widetilde{\text{Id}}_i$ be the
  corresponding lift of the identity map $\text{Id}$ to $\tilde V_i$
  and let $B^n(r)$ denote the metric $r$-ball centered at $0$ in
  $\R^n$.  Choose $\varepsilon_i$ small enough so that if $\tilde
  D_i=\widetilde{\text{Id}}_i^{-1}(B(1-\varepsilon_i))$ then the
  collection $\{D_i=\pi_i(\tilde D_i)\}$ covers $\orbify{O}$ and also
  that $B(1+\varepsilon_i)\subset\tilde V_i$.

  Let $f:\orbify{O}\to\orbify{O}$ be a map such that $\|\tilde
  f_i(\tilde x)-\widetilde{\text{Id}}_i(\tilde x)\|_{\tilde
    V_i}<\varepsilon_i$ for $\tilde x\in\tilde C_i$ and all $i$. We
  want to show that $f$ is surjective.

  Define $\tilde g_i=\tilde f_i\circ\widetilde{\text{Id}}_i^{-1}$.
  Then $\tilde g_i$ is a map from $B^n=\tilde C_i$ into $\R^n$ and the
  image of the unit sphere $S^{n-1}=\partial B^n$ under $\tilde g_i$
  lies outside $B(1-\varepsilon_i)$. We will show that $\tilde
  D_i\subset\tilde g_i(B^n)$.  Since $\{D_i\}$ cover $\orbify{O}$ and
  $D_i=\pi_i(\tilde D_i)\subset\pi_i\circ\tilde
  g_i(B^n)=\pi_i\circ\tilde f_i(\tilde C_i)=f(C_i)$, this will imply
  that $f$ is surjective.

  Suppose to the contrary that $\tilde y\in B(1-\varepsilon_i)$, but
  $\tilde y\notin\tilde g_i(B^n)$. Let $\lambda:\R^n-\{\tilde y\}\to
  S^{n-1}$ be the radial projection from $\tilde y$. Then
  $\lambda\circ\tilde g_i$ maps $B^n$ into $S^{n-1}$. On the other
  hand, the restriction $\tilde g_i|_{S^{n-1}}:S^{n-1}\to\R^n$ is
  homotopic to the identity map via $F_t(\tilde x)=t\tilde g_i(\tilde
  x)+(1-t)\tilde x$ for $\tilde x\in S^{n-1}$. This homotopy carries
  $\tilde g_i(\tilde x)$ along the straight line between $\tilde
  g_i(\tilde x)$ and $\tilde x$ so $F_t(\tilde x)$ lies outside
  $B(1-\varepsilon_i)$. Thus, $\lambda\circ F_t$ is a well-defined
  homotopy between $(\lambda\circ\tilde g_i)|_{S^{n-1}}:S^{n-1}\to
  S^{n-1}$ and the identity map. It is not necessary that $F_t$ and
  $\lambda$ be equivariant.  Now consider the homology sequence of the
  pair $(B^n,S^{n-1})$:

\begin{equation*}
  \xymatrix{
    {0}\ar[r]&{H_n(B^n,S^{n-1})}\ar[r]\ar[d]_{(\lambda\circ\tilde g_i)_*}&%
    {H_{n-1}(S^{n-1})}\ar[r]\ar[d]^{\left((\lambda\circ\tilde g_i)|_{S^{n-1}}\right)_*} & {0} \\
    {0}\ar[r] & {H_n(B^n,S^{n-1})}\ar[r] & {H_{n-1}(S^{n-1})}\ar[r] & {0}  }
\end{equation*}
$(\lambda\circ\tilde g_i)_*$ is the zero homomorphism since
$(\lambda\circ\tilde g_i)$ sends $B^n$ into $S^{n-1}$. However,
$\left((\lambda\circ\tilde g_i)|_{S^{n-1}}\right)_*$ is the identity
homomorphism since $(\lambda\circ\tilde g_i)|_{S^{n-1}}$ is homotopic
to the identity map. Since $H_n(B^n,S^{n-1})\cong\Z$ and the diagram
commutes we have a contradiction. Thus, $f$ is surjective.

\end{proof}

The following is a culmination of the results of this section.

\begin{theorem}\label{ExpSmallOrbiIsDiffeo} Let $\orbify{O}$ be a
  $C^{r+1}$ Riemannian orbifold.  If $\sigma$ is a $C^r$ orbisection
  sufficiently $C^1$ close to the $0$-orbisection $\mathbf{0}$ of
  $T\orbify{O}$ then $E^\sigma$ is a $C^r$ orbifold diffeomorphism.
  That is, $E^\sigma\in\Diff^r_{\Orb}(\orbify{O})$.
\end{theorem}
\begin{proof} Let $\{C_i\}$ be a locally finite covering of
  $\orbify{O}$ by compact sets. By
  proposition~\ref{ExpOrbisectionIsLocalDiffeo}, there exist positive
  constants $\varepsilon_i$ such that if $\sigma$ is $C^1$
  $\varepsilon_i$-close to $\mathbf{0}$ on $C_i$, then
  $E^\sigma|_{C_i}$ is a $C^r$ orbifold embedding. Since
  $\text{Id}=E^{\mathbf{0}}=(\exp\circ\mathbf{0})$, by choosing
  $\varepsilon_i$ smaller if necessary, we may conclude that
  $E^\sigma$ is surjective by lemma~\ref{C0CloseToIdIsOnto}. We need
  only to show that $E^\sigma$ is globally injective. To do this, we
  modify the argument in \cite[theorem 3.10]{MR0198479}.

  Let $\{D_i\}$ be a covering of $\orbify{O}$ by compact sets with
  $D_i\subset\intr(C_i)$. Let
  $\delta_i=d\left(D_i,\orbify{O}-\intr(C_i)\right)>0$.  By choosing
  $\varepsilon_i$ smaller if necessary, we may assume that $E^\sigma$
  is $C^1$ $\frac12\delta_i$-close to $\text{Id}$ for $x\in D_i$ and
  that $E^\sigma(D_i)\subset C_i$. Suppose that
  $E^\sigma(x)=E^\sigma(y)$, where $x\in D_i$ and $y\in D_j$ and
  $\delta_i\le\delta_j$. Then
$$d(x,y)\le d(x, E^\sigma(x))+d(E^\sigma(x),E^\sigma(y))+d(E^\sigma(y),y)<\frac12\delta_i+\frac12\delta_j\le\delta_j$$
However, since $E^\sigma$ is injective on $C_j$, $x\notin C_j$. Thus,
$d(x,y)\ge\delta_j$, a contradiction. Hence $E^\sigma$ is injective
and thus a $C^r$ orbifold diffeomorphism.
\end{proof}

\section{Proof of Theorem~1 and Corollary~2}\label{ProofOfMainTheoremSection}

Throughout this section, we assume that $\orbify{O}$ is a smooth
compact orbifold (without boundary). Without loss of generality, we
may assume, by propositions~\ref{OrbifoldsAreSmoothable} and
\ref{OrbifoldRiemannianMetricExists}, that $\orbify{O}$ is a
$C^\infty$ orbifold with $C^\infty$ Riemannian metric. We let
$\mathscr{B}^r(\sigma,\varepsilon)=\mathscr{N}^r(\sigma,\varepsilon)\cap\mathscr{D}^r(\orbify{O})$.
That is, $\mathscr{B}^r(\sigma,\varepsilon)$ is the set of $C^r$
orbisections $\varepsilon$-close to $\sigma$ in the $C^r$ topology on
$C^r_{\Orb}(\orbify{O},T\orbify{O})$. We prove the main theorem in a
series of propositions.

\begin{proposition} There exists $\varepsilon>0$ such that
  $E^\sigma=\exp\circ\sigma\in\Diff^r_{\Orb}(\orbify{O})$ for
  $\sigma\in\mathscr{B}^r(\mathbf{0},\varepsilon)$. That is, there
  exists a map
  $E:\mathscr{B}^r(\mathbf{0},\varepsilon)\to\Diff^r_{\Orb}(\orbify{O})$
  defined by $E(\sigma)=E^\sigma$.
\end{proposition}
\begin{proof} This follows from compactness of $\orbify{O}$ and
  theorem~\ref{ExpSmallOrbiIsDiffeo}.
\end{proof}

\begin{proposition}\label{EIsInjective} The map
  $E:\mathscr{B}^r(\mathbf{0},\varepsilon)\to\Diff^r_{\Orb}(\orbify{O})$
  is injective.
\end{proposition}
\begin{proof} Suppose $E(\sigma)=E(\tau)$ for
  $\sigma,\tau\in\mathscr{B}^r(\mathbf{0},\varepsilon)$. Then
  $(\exp\circ\sigma)(x)=(\exp\circ\tau)(x)$ for all $x\in\orbify{O}$.
  Thus, in each orbifold chart $(\tilde U_x,\Gamma_x)$, we have
  $\pi_x\circ\widetilde{\exp}_{\tilde U_x}(\tilde x,\tilde
  v)=\pi_x\circ\widetilde{\exp}_{\tilde U_x}(\tilde x,\tilde w)$.
  Since $\widetilde{\exp}_{\tilde U_x}(\tilde x,\cdot)$ is a local
  $C^r$ diffeomorphism we must have $\tilde v=(d\gamma)_{\tilde
    x}(\tilde w)$ for some $\gamma\in\Gamma_x$. Thus, $v=w\in
  A_x\orbify{O}$. Hence $\sigma=\tau$ and $E$ is injective.
\end{proof}

\begin{proposition}\label{EIsSurjective} The map 
  $E:\mathscr{B}^r(\mathbf{0},\varepsilon)\to\mathscr{N}^0(\textup{Id},\varepsilon)\cap\Diff^r_{\Orb}(\orbify{O})$
  is surjective.
\end{proposition}
\begin{proof}
  Let
  $f\in\mathscr{N}^0(\textup{Id},\varepsilon)\cap\Diff^r_{\Orb}(\orbify{O})$.
  Let $\{C_i\}$ be a finite covering of $\orbify{O}$ by compact sets
  such that $C_i$ is an orbifold chart and $f(C_i)\subset V_i$ where
  $V_i$ is a relatively compact orbifold chart. Let $x\in C_i$, and
  $\tilde U_x\subset\intr{\tilde C_i}$ an orbifold chart at $x$ where
  the local lift $\tilde f_x$ to $\tilde U_x$ is $C^0$
  $\varepsilon$-close to the lift
  $\widetilde{\text{Id}}_x=\text{Id}_{\tilde U_x}$ of the identity map
  and not $\varepsilon$-close to any other lift of the identity over
  $\tilde U_x$.  For $\varepsilon$ small enough it follows that
  $\Theta_{f,x}(\delta)=\Theta_{\textup{Id},x}(\delta)=\delta$ for all
  $\delta\in\Gamma_x$. This is because for each $\delta\in\Gamma_x$ we
  have

\begin{align*}
  \|\tilde f_x(\delta\cdot\tilde y)-\widetilde{\text{Id}}_x(\delta\cdot\tilde y)\|_{\tilde V_i}<\varepsilon &\Longleftrightarrow\\
  \|\Theta_{f,x}(\delta)\cdot\tilde f_x(\tilde y)-\delta\cdot\tilde y\|_{\tilde V_i}<\varepsilon & \Longleftrightarrow\\
  \|\delta^{-1}\Theta_{f,x}(\delta)\cdot\tilde f_x(\tilde y)-\tilde
  y\|_{\tilde V_i}<\varepsilon %
  & \Longleftrightarrow\text{(since }\Gamma_x\text{ acts by isometries)}\\
  \|\delta^{-1}\Theta_{f,x}(\delta)\cdot\tilde f_x(\tilde
  y)-\widetilde{\text{Id}}_x(\tilde y)\|_{\tilde V_i}<\varepsilon
\end{align*}
Thus, by our choice of local lift of the identity map over $\tilde
U_x$, it follows that $\delta^{-1}\Theta_{f,x}(\delta)=e$ which
implies that $\Theta_{f,x}(\delta)=\delta$.

We wish to define a $C^r$ orbisection $\sigma$ so that $E(\sigma)=f$.
We do this by defining appropriate local lifts $\tilde\sigma_x$. In
particular, let
$$\tilde\sigma_x(\tilde y)=%
\left(\tilde y,\widetilde{\exp}^{-1}_{\tilde U_x,\tilde y}\left(\tilde
    f_x(\tilde y)\right)\right)\in T\tilde U_x$$ Before we show that
$\tilde\sigma_x$ satisfies the correct equivariance relation observe
that, in general, $\widetilde{\exp}^{-1}_{\tilde U_x,\tilde
  y}\left(\gamma\cdot\tilde z\right)=%
(d\gamma)_{\gamma^{-1}\tilde y}\circ\widetilde{\exp}^{-1}_{\tilde
  U_x,\gamma^{-1}\tilde y}(\tilde z)=%
\gamma\cdot\widetilde{\exp}^{-1}_{\tilde U_x,\gamma^{-1}\tilde
  y}(\tilde z)$. Thus,

\begin{align*}
  \tilde\sigma_x(\delta\cdot\tilde y) & =
  \left(\delta\cdot\tilde y,\widetilde{\exp}^{-1}_{\tilde
      U_x,\delta\tilde y}%
    \left(\tilde f_x(\delta\cdot\tilde y)\right)\right)\\
  & = \left(\delta\cdot\tilde y,\widetilde{\exp}^{-1}_{\tilde
      U_x,\delta\tilde y}%
    \left(\delta\cdot\tilde f_x(\tilde y)\right)\right)\\
  & = \left(\delta\cdot\tilde y,\delta\cdot%
    \widetilde{\exp}^{-1}_{\tilde U_x,\delta^{-1}\delta\tilde y}%
    \left(\tilde f_x(\tilde y)\right)\right)\\
  & = \left(\delta\cdot\tilde y,\delta\cdot%
    \widetilde{\exp}^{-1}_{\tilde U_x,\tilde y}\left(\tilde f_x(\tilde y)\right)\right)\\
  & = \delta\cdot\tilde\sigma_x(\tilde y)
\end{align*}
which is the correct equivariance relation for an orbisection. As a
result we see that the map $\sigma(x)=\Pi_x\circ\tilde\sigma_x(\tilde
x)$ defines a $C^r$ orbisection of $T\orbify{O}$ and that
$E(\sigma)=f$ since $\tilde\sigma_x(\tilde x)=%
\left(\tilde x,\widetilde{\exp}^{-1}_{\tilde U_x,\tilde x}\left(\tilde
    f_x(\tilde x)\right)\right)$.
\end{proof}

The following proposition is the last ingredient needed to complete the proof
of theorem~\ref{MainTheorem}.

\begin{proposition}\label{EIsAHomeomorphism} The map
  $E:\mathscr{B}^r(\mathbf{0},\varepsilon)\to\mathscr{N}^0(\textup{Id},\varepsilon)\cap\Diff^r_{\Orb}(\orbify{O})$
  is a homeomorphism.
\end{proposition}
\begin{proof} Propositions ~\ref{EIsInjective} and \ref{EIsSurjective}
  show that $E$ is bijective. Continuity of $E$ follows from the
  formula for a local lift of $E$ given in
  Propositon~\ref{ExpOrbisectionIsMap} and continuity of $E^{-1}$
  follows from the formula for $\tilde\sigma_x$ given in the last line
  of proposition~\ref{EIsSurjective}.
\end{proof}

\subsection*{\bf Proof of Theorem \ref{MainTheorem}}

\begin{proof} Let $f\in\Orbdiff^r(O)$. By proposition~\ref{EIsAHomeomorphism}, the map
$$f\circ E:\mathscr{B}^r(\mathbf{0},\varepsilon)\to\mathscr{N}^0(f,\varepsilon)\cap\Diff^r_{\Orb}(\orbify{O})$$
is a homeomorphism giving a local chart about $f$. Let
$\mathscr{N}_{fg}=\mathscr{N}^0(f,\varepsilon)\cap\mathscr{N}^0(g,\varepsilon)\cap\Diff^r_{\Orb}(\orbify{O})$ denote a chart overlap, and let
$\mathscr{B}_{fg}=(f\circ E)^{-1}(\mathscr{N}_{fg})\subset\mathscr{B}^r(\mathbf{0},\varepsilon)$. Then
the corresponding transition map
$$\left.(g\circ E)^{-1}\circ (f\circ E)\right|_{\mathscr{B}_{fg}}: %
\mathscr{B}_{fg}\subset\mathscr{B}^r(\mathbf{0},\varepsilon)\to %
(g\circ E)^{-1}(\mathscr{N}_{fg})\subset\mathscr{B}^r(\mathbf{0},\varepsilon)$$
is a homeomorphism. This gives the desired $C^0$ manifold structure to
$\Diff^r_{\Orb}(\orbify{O})$ where the model space is the topological
vector space of $C^r$ orbisections of the tangent orbibundle with the
$C^r$ topology.
\end{proof}

\subsection*{Proof of Corollary \ref{MainCorollary}}

\begin{proof}  It follows from the arguments in example \ref{IdentityMap} that for
  a given $f\in \mathscr{ID}$ and any $x\in\orbify{O}$ with orbifold
  chart $U_x$ of $x$ there is a $\gamma_x\in\Gamma_x$ so that $\tilde
  f(\tilde y) = \gamma_x\cdot\tilde y$ for all $\tilde y\in\tilde
  U_x$.  A finite cover of $\orbify{O}$ by charts
  $\{U_{x_1},\ldots,U_{x_M}\}$ shows that $\mathscr{ID}$ is a subgroup
  of $\prod_{i=1}^M\Gamma_{x_i}$ and is therefore finite.  Clearly
  $\mathscr{ID}$ is a normal subgroup of $\Diff^r_{\Orb}(\orbify{O})$
  as $\tilde g\circ\tilde f\circ\tilde g^{-1}$ covers the identity for
  any $g\in\Diff^r_{\Orb}(\orbify{O})$ and $f\in \mathscr{ID}$. Also, any
  two lifts $\tilde h_0$ and $\tilde h_1$ of
  $h\in\Diff^r_{\textup{red}}(\orbify{O})$ by definition must satisfy
  $\tilde h_0\circ\tilde h_1^{-1}\in \mathscr{ID}$ from which follows
  the existence of the short exact sequence.  Moreover, the finiteness
  of $\mathscr{ID}$ shows that the quotient topology on
  $\Diff^r_{\textup{red}}(\orbify{O})$ is again that of a Banach
  manifold if $r < \infty$ and of a \Frechet manifold if $r = \infty$.
\end{proof}

\bibliography{ref}
\bibliographystyle{amsalpha}

\end{document}